\pdfoutput=1

\documentclass[11pt]{amsart}
\usepackage[hmargin=3cm,vmargin=3cm]{geometry}

\usepackage[latin1]{inputenc}
\usepackage{xspace,amssymb,amsfonts,euscript}
\usepackage{amsthm,amsmath,amscd,stmaryrd,latexsym}
\usepackage{palatino, euscript}

\usepackage{graphicx}
\usepackage[all]{xy}
\usepackage[dvips]{epsfig}
\usepackage{pinlabel}
\usepackage{psfrag}

\usepackage{float}
\usepackage{blkarray}

\usepackage{tikz}
\usetikzlibrary{cd}
\usetikzlibrary{decorations.markings}
\usetikzlibrary{decorations.pathreplacing}
\usetikzlibrary{arrows,shapes,positioning}
\usetikzlibrary{decorations.pathmorphing}
\tikzset{anchorbase/.style={baseline={([yshift=-0.5ex]current bounding box.center)}},
tinynodes/.style={font=\tiny,text height=0.75ex,text depth=0.15ex},
smallnodes/.style={font=\scriptsize,text height=0.75ex,text depth=0.15ex},
>={Latex[length=1mm, width=1.5mm]},
overcross/.style={line width=4pt,color=white},
}
\tikzstyle directed=[postaction={decorate,decoration={markings,
    mark=at position #1 with {\arrow{>}}}}]
\tikzstyle rdirected=[postaction={decorate,decoration={markings,
    mark=at position #1 with {\arrow{<}}}}]
    
\IfFileExists{comments.sty}{\usepackage{comments}}

\RequirePackage{color}
\definecolor{myred}{rgb}{0.75,0,0}
\definecolor{mygreen}{rgb}{0,0.5,0}
\definecolor{myblue}{rgb}{0,0,0.65}

\RequirePackage{ifpdf}
\ifpdf
  \IfFileExists{pdfsync.sty}{\RequirePackage{pdfsync}}{}
  \RequirePackage[pdftex,
   colorlinks = true,
   urlcolor = myblue, 
   citecolor = mygreen,  
   linkcolor = myred, 
   pagebackref,
   bookmarksopen=true]{hyperref}
\else
  \RequirePackage[hypertex]{hyperref}
\fi

\newtheorem{thm}{Theorem}[section]
\newtheorem{lemma}[thm]{Lemma}

\newtheorem{prop}[thm]{Proposition}
\newtheorem{cor}[thm]{Corollary}
  
\newtheorem{conj}[thm]{Conjecture}

\newtheorem*{prop*}{Proposition}

\theoremstyle{definition}
\newtheorem{defn}[thm]{Definition}

\newtheorem{notation}[thm]{Notation}

\newtheorem{example}[thm]{Example}

\theoremstyle{remark}
\newtheorem{remark}[thm]{Remark}
\newtheorem{rmk}[thm]{Remark}

\numberwithin{equation}{section}

     \def\ak{{\mathbbm{k}}}

    \def\DM{{\mathbb{D}}}

\let\phi=\varphi

\usepackage{bbm}

\def\1{\mathbbm{1}}

\newcommand{\un}{\underline}

\newcommand{\ot}{\otimes}

\renewcommand{\to}{\rightarrow}

\newcommand{\refequal}[1]{\xy {\ar@{=}^{#1}
(-1,0)*{};(1,0)*{}};
\endxy}

\DeclareMathOperator{\Kar}{\bf{Kar}}
\DeclareMathOperator{\Hom}{{\rm Hom}}

\DeclareMathOperator{\End}{{\rm End}}

\DeclareMathOperator{\Ext}{{\rm Ext}}

\DeclareMathOperator{\id}{{\rm id}}
\DeclareMathOperator{\im}{{\rm im}}

\DeclareMathOperator{\DD}{\mathcal{D}_{\mathfrak{sp}_4}}
\DeclareMathOperator{\DDk}{\mathcal{D}^{\ak}_{\mathfrak{sp}_4}}
\DeclareMathOperator{\eval}{\Xi}

\DeclareMathOperator{\Fund}{\bf{Fund}}
\DeclareMathOperator{\Rep}{\bf{Rep}}

\DeclareMathOperator{\Tilt}{\bf{Tilt}}
\DeclareMathOperator{\wt}{\rm{wt}}

\newcommand{\blues}{\textcolor{blue}{1}}
\newcommand{\greent}{\textcolor{green}{2}}

\begin{document}

\begin{abstract} Using the light ladder basis for Kuperberg's $C_2$ webs, we derive triple clasp formulas for idempotents projecting to the top summand in each tensor product of fundamental representations. We then find explicit formulas for the coefficients occurring in the clasps, by computing these coefficients as local intersection forms. Our formulas provide further evidence for Elias's clasp conjecture, which was given for type $A$ webs, and suggests how to generalize the conjecture to non-simply laced types. 
\end{abstract}

\title{Triple Clasp Formulas for $C_2$ Webs}

\author{Elijah Bodish} \address{University of Oregon, Eugene}

\maketitle



\section{Introduction}
\label{sec-intro}

Let $\mathfrak{g}$ be a semisimple Lie algebra. Fix a dominant integral weight $\lambda$. We can uniquely write $\lambda = \sum_{\varpi} \lambda_{\varpi}\varpi$, where the $\varpi$ are fundamental weights. The irreducible representation with highest weight $\lambda$, denoted $V(\lambda)$, occurs with multiplicity one as a direct summand in the tensor product 
\begin{equation}\label{topsumnd}
\bigotimes_{\varpi}V(\varpi)^{\otimes \lambda_{\varpi}},
\end{equation}
and all other irreducible summands are isomorphic to $V(\mu)$, with $\mu < \lambda$ in the dominance order. Moreover, the isomorphism class of the tensor product is unaffected by the order of the tensor factors in Equation \eqref{topsumnd}. For any sequence of fundamental weights: $\varpi_1, \varpi_2, \ldots, \varpi_d$, we refer to
\begin{equation}
V(\varpi_1 + \varpi_2 + \ldots + \varpi_d)\subset V(\varpi_1)\otimes V(\varpi_2)\otimes \ldots \otimes V(\varpi_d)
\end{equation}
as the \emph{top summand}. 

Our main result is to use Kuperberg's type $B_2= C_2$ webs to give a recursive description of the idempotent projecting to the top summand in an arbitrary tensor product of fundamental representations of $\mathfrak{sp}_4$. Our formulas are the $\mathfrak{sp}_4$ analogue of Elias's conjectural recursive formulas describing the idempotent projecting to the top summand for $\mathfrak{sl}_n$ \cite{elias2015light}. 

Elias's formulas were inspired by Wenzl's recursive formula for $\mathfrak{sl}_2$. The $\mathfrak{sl}_2$ case so well illustrates the arguments used to derive our main Theorem \eqref{MAINTHM} that we recall Wenzl's recursion below. Experts can skip to Section \eqref{sec-results}. 

\subsection{Wenzl's triple clasp formula for $A_1$ webs}
\label{sec-sl2}

Let $\mathbb{C}^2$ be the two dimensional defining representation of $\mathfrak{sl}_2(\mathbb{C})$. The tensor powers $(\mathbb{C}^2)^{\otimes d}$ carry an action of $\mathfrak{sl}_2(\mathbb{C})$ and the module $S^d(\mathbb{C}^2)$ is an irreducible quotient of $(\mathbb{C}^2)^{\otimes d}$. The finite dimensional irreducible representations of $\mathfrak{sl}_2$ are in bijection with $\mathbb{Z}_{\ge 0}$ via $d\mapsto S^d(\mathbb{C}^2)$ and the composition factors of the kernel of the map $(\mathbb{C}^2)^{\otimes d}\longrightarrow S^d(\mathbb{C}^2)$ are all of the form $S^k(\mathbb{C}^2)$ for $k< d$. 

Let $\Rep(\mathfrak{sl}_2(\mathbb{C}))$ denote the abelian monoidal category of finite dimensional representations of $\mathfrak{sl}_2(\mathbb{C})$. Since $\Rep(\mathfrak{sl}_2(\mathbb{C}))$ is semisimple and has simple objects in bijection with the nonnegative integers, it is equivalent to $\oplus_{\mathbb{Z}_{\ge 0}} \mathbf{Vec}_{\mathbb{C}}$, and is therefore uninteresting as an abelian category. However, a semisimple monoidal category contains much more information than just the number of simple objects. For example, by the general theory of Tannakian reconstruction \cite{delignemilne} one can recover the group $SL_2(\mathbb{C})$ as the automorphisms of the monoidal functor $\Fund(\mathfrak{sl}_2(\mathbb{C}))\longrightarrow \mathbf{Vec}_{\mathbb{C}}$.

Define $\Fund(\mathfrak{sl}_2(\mathbb{C}))$ to be the full monoidal subcategory of $\Rep(\mathfrak{sl}_2(\mathbb{C}))$ with objects arbitrary tensor products of $\mathbb{C}^2$. Since each irreducible finite dimensional representation of $\mathfrak{sl}_2(\mathbb{C})$ is a direct summand of a tensor power of the defining representation, there is an equivalence of monoidal categories:
 \[
 \Kar(\Fund(\mathfrak{sl}_2(\mathbb{C}))) \cong \Rep(\mathfrak{sl}_2(\mathbb{C})).\footnote{The Karoubi envelope of a category $\mathcal{C}$, denoted $\Kar(\mathcal{C})$ is the category with objects pairs: $(X, e)$, where $X$ is an object in $\mathcal{C}$ and $e\in \End_{\mathcal{C}}(X)$ is an idempotent. The morphisms $(X, e)\longrightarrow (Y, f)$ in $\Kar(\mathcal{C})$ are all morphisms of the form $f\circ \varphi\circ e$, where $\varphi:X\longrightarrow Y$ in $\mathcal{C}$. When $\mathcal{C}$ is a monoidal category, $\Kar(\mathcal{C})$ is also naturally a monoidal category. Moreover, if $\mathcal{A}$ is a semisimple monoidal category, $\mathcal{C}\subset \mathcal{A}$ is a full monoidal subcategory, and every irreducible object in $\mathcal{A}$ is a direct summand of an object in $\mathcal{C}$, then there is a monoidal equivalence $\Kar(\mathcal{C})\longrightarrow \mathcal{A}$.} 
 \]
Threfore, the study of $\Rep(\mathfrak{sl}_2(\mathbb{C}))$ as a monoidal category is reduced to the study of idempotents in $\Fund(\mathfrak{sl}_2(\mathbb{C}))$. 

Let $\mathcal{TL}$ be the strict, pivotal, and $\mathbb{C}$-linear category generated by one self dual object of dimension $-2$. It is well known that $\mathcal{TL}$ is equivalent to the monoidal category $\Fund(\mathfrak{sl}_2(\mathbb{C}))$ \cite{RTW, TemLie}. Thus, we are led to consider the problem of using the category $\mathcal{TL}$ to describe the idempotent in $\End_{\mathfrak{sl}_2(\mathbb{C})}((\mathbb{C}^2)^{\ot d})$ which projects to $S^d(\mathbb{C}^2)$. Wenzl found a recursive description of these idempotents \cite{wenzl_1987}. Using the usual graphical calculus for $\mathcal{TL}$ and using a $d$ labelled oval to represent the idempotent with image $S^d(\mathbb{C}^2)$, the Wenzl recursion becomes the following. 
\begin{equation}\label{JW}
\begin{tikzpicture}[scale=.4, anchorbase]
	\draw[ultra thick, black] (-2,-3.2) -- (-2,3.2);
	\draw[ultra thick, black] (2, -3.2) -- (2, 3.2);
	\draw[ultra thick, black] (1,-3.2) -- (1, 3.2);
	\draw[ultra thick, black, fill= white] (0,0) ellipse (2.5 and 1) node{$d+1$};
\end{tikzpicture}
\quad = \quad
\begin{tikzpicture}[scale=.4, anchorbase]
	\draw[ultra thick, black] (-2,-3.2) -- (-2,3.2);
	\draw[ultra thick, black] (2, -3.2) -- (2, 3.2);
	\draw[ultra thick, black] (1,-3.2) -- (1, 3.2);
	\draw[ultra thick, black, fill= white] (-.55,0) ellipse (2 and .7) node{$d$};
\end{tikzpicture}
\quad + \dfrac{d}{d+1}\quad
\begin{tikzpicture}[scale=.4, anchorbase]
	\draw[ultra thick, black] (-2,-3.2) -- (-2,3.2);
	\draw[ultra thick, black] (2, -3.2) -- (2, -1.4);
	\draw[ultra thick, black] (2, 1.4) -- (2, 3.2);
	\draw[ultra thick, black] (1,-3.2) -- (1, -1.4);
	\draw[ultra thick, black] (1, 1.4) -- (1, 3.2);
	\draw[ultra thick, black] (2,1.4) arc (0: -180: .5);
	\draw[ultra thick, black] (2,-1.4) arc (0: 180: .5);
	\draw[ultra thick, black, fill= white] (-.55,-2.2) ellipse (2 and .7) node{$d$};
	\draw[ultra thick, black, fill= white] (-.55,2.2) ellipse (2 and .7) node{$d$};
	\draw[ultra thick, black, fill= white] (0,0) ellipse (2.3 and .7) node{$d-1$};
\end{tikzpicture}
\end{equation}

One way to think about Wenzl's recursion, which we learned from \cite{elias2015light}, is to first note that by Schur's lemma there is some $\kappa_d\in \mathbb{C}$ such that the following equality holds.
\begin{equation}\label{sl2lif}
\begin{tikzpicture}[scale=.4, anchorbase]
	\draw[ultra thick, black] (-2,-3.2) -- (-2,3.2);
	\draw[ultra thick, black] (1, 0) circle (1);
	\draw[ultra thick, black, fill= white] (-1,0) ellipse (2 and .7) node{$d$};
	\draw[ultra thick, black, fill= white] (0,3) ellipse (2.3 and .7) node{$d-1$};
	\draw[ultra thick, black, fill= white] (0,-3) ellipse (2.3 and .7) node{$d-1$};
\end{tikzpicture}
\quad = \kappa_d \quad
\begin{tikzpicture}[scale=.4, anchorbase]
	\draw[ultra thick, black] (-2,-3.2) -- (-2,3.2);
	\draw[ultra thick, black, fill= white] (0,0) ellipse (2.3 and .7) node{$d-1$};
\end{tikzpicture}
\end{equation}
Since
\[
S^d(\mathbb{C}^2)\otimes \mathbb{C}^2\cong S^{d+1}(\mathbb{C}^2)\oplus S^{d-1}(\mathbb{C}^2),
\]
we also observe that there is a relation in $\End_{\mathfrak{sl}_2(\mathbb{C}))}(\mathbb{C}^2)$ of the following form. 
\begin{equation}\label{sl2recursion}
\begin{tikzpicture}[scale=.4, anchorbase]
	\draw[ultra thick, black] (-2,-3.2) -- (-2,3.2);
	\draw[ultra thick, black] (2, -3.2) -- (2, 3.2);
	\draw[ultra thick, black] (1,-3.2) -- (1, 3.2);
	\draw[ultra thick, black, fill= white] (-.55,0) ellipse (2 and .7) node{$d$};
\end{tikzpicture}
\quad = \quad
\begin{tikzpicture}[scale=.4, anchorbase]
	\draw[ultra thick, black] (-2,-3.2) -- (-2,3.2);
	\draw[ultra thick, black] (2, -3.2) -- (2, 3.2);
	\draw[ultra thick, black] (1,-3.2) -- (1, 3.2);
	\draw[ultra thick, black, fill= white] (0,0) ellipse (2.5 and 1) node{$d+1$};
\end{tikzpicture}
\quad + \kappa_d^{-1}\quad
\begin{tikzpicture}[scale=.4, anchorbase]
	\draw[ultra thick, black] (-2,-3.2) -- (-2,3.2);
	\draw[ultra thick, black] (2, -3.2) -- (2, -1.4);
	\draw[ultra thick, black] (2, 1.4) -- (2, 3.2);
	\draw[ultra thick, black] (1,-3.2) -- (1, -1.4);
	\draw[ultra thick, black] (1, 1.4) -- (1, 3.2);
	\draw[ultra thick, black] (2,1.4) arc (0: -180: .5);
	\draw[ultra thick, black] (2,-1.4) arc (0: 180: .5);
	\draw[ultra thick, black, fill= white] (-.55,-2.2) ellipse (2 and .7) node{$d$};
	\draw[ultra thick, black, fill= white] (-.55,2.2) ellipse (2 and .7) node{$d$};
	\draw[ultra thick, black, fill= white] (0,0) ellipse (2.3 and .7) node{$d-1$};
\end{tikzpicture}
\end{equation}
Where $\kappa_d^{-1}$ is the coefficient needed to make the quasi-idempotent 
\[
\begin{tikzpicture}[scale=.4, anchorbase]
	\draw[ultra thick, black] (-2,-3.2) -- (-2,3.2);
	\draw[ultra thick, black] (2, -3.2) -- (2, -1.4);
	\draw[ultra thick, black] (2, 1.4) -- (2, 3.2);
	\draw[ultra thick, black] (1,-3.2) -- (1, -1.4);
	\draw[ultra thick, black] (1, 1.4) -- (1, 3.2);
	\draw[ultra thick, black] (2,1.4) arc (0: -180: .5);
	\draw[ultra thick, black] (2,-1.4) arc (0: 180: .5);
	\draw[ultra thick, black, fill= white] (-.55,-2.2) ellipse (2 and .7) node{$d$};
	\draw[ultra thick, black, fill= white] (-.55,2.2) ellipse (2 and .7) node{$d$};
	\draw[ultra thick, black, fill= white] (0,0) ellipse (2.3 and .7) node{$d-1$};
\end{tikzpicture}
\]
into an idempotent.
Using the $d\mapsto d-1$ version of the relation in Equation \eqref{sl2recursion} to rewrite the middle clasp labelled $d$ on the left hand side of Equation \eqref{sl2lif}, we find the following. 
\begin{equation}\label{resolvesl2}
\begin{tikzpicture}[scale=.4, anchorbase]
	\draw[ultra thick, black] (-2,-3.2) -- (-2,3.2);
	\draw[ultra thick, black] (2.5, 0) circle (1);
	\draw[ultra thick, black, fill= white] (0,0) ellipse (2.5 and 1) node{$d$};
\end{tikzpicture}
\quad = \quad
\begin{tikzpicture}[scale=.4, anchorbase]
	\draw[ultra thick, black] (-2,-3.2) -- (-2,3.2);
	\draw[ultra thick, black] (3, 0) circle (1);
	\draw[ultra thick, black, fill= white] (-.55,0) ellipse (2 and .7) node{$d-1$};
\end{tikzpicture}
\quad - \kappa_{d-1}^{-1}\quad
\begin{tikzpicture}[scale=.4, anchorbase]
	\draw[ultra thick, black] (-2,-3.2) -- (-2,3.2);
	\draw[ultra thick, black] (2, -3.2) -- (2, -1.4);
	\draw[ultra thick, black] (3, -3.2) -- (3, 3.2);
	\draw[ultra thick, black] (2, 1.4) -- (2, 3.2);
	\draw[ultra thick, black] (1,-3.2) -- (1, -1.4);
	\draw[ultra thick, black] (1, 1.4) -- (1, 3.2);
	\draw[ultra thick, black] (2,1.4) arc (0: -180: .5);
	\draw[ultra thick, black] (2,-1.4) arc (0: 180: .5);
	\draw[ultra thick, black] (2,-3.2) arc (-180: 0: .5);
	\draw[ultra thick, black] (2,3.2) arc (180: 0: .5);
	\draw[ultra thick, black, fill= white] (-.55,-2.2) ellipse (2 and .7) node{$d-1$};
	\draw[ultra thick, black, fill= white] (-.55,2.2) ellipse (2 and .7) node{$d-1$};
	\draw[ultra thick, black, fill= white] (0,0) ellipse (2.3 and .7) node{$d-2$};
\end{tikzpicture}
\end{equation}
Using the relations in $\mathcal{TL}$, along with the fact that post-composing the ovals with any cap map results in zero, we can simplify the right hand side of Equation \eqref{resolvesl2} and find $\kappa_d= -2 - 1/\kappa_{d-1}$. From the initial condition $\kappa_0^{-1} = 0$, it is easy to verify that $\kappa_d= -(d+1)/d$.

\begin{notation}
The terminology of \emph{clasp} was introduced in \cite{Kupe} to refer to idempotents projecting to the top summand expressed in terms of the graphical calculus. Following \cite{elias2015light}, we will refer to a recursive formula for clasps, in which the terms on the right hand side of the recursion are three clasps linked together by diagrammatic morphisms, as a \emph{triple clasp formula}. We will refer to the diagrammatic morphisms in the clasp formula as \emph{elementary light ladder diagrams}.
\end{notation}

\begin{remark}
Equation \eqref{JW} is still true without the middle clasp labelled $d-1$. However, the middle clasp in the triple clasp keeps track of which summand the morphism is factoring through and therefore has representation theoretic meaning. Moreover, in higher rank examples, like the one considered in this article, removing the middle clasp will not result in a valid identity. One would have to also change the coefficients and we do not expect these new coefficients to be as nice as those occurring in the triple clasp formula.
\end{remark}

\begin{remark}
The category $\Rep(\mathfrak{sl}_2(\mathbb{C}))$ has a $\mathbb{C}(q)$-linear analogue, the category of finite dimensional type $1$ representations of $U_q(\mathfrak{sl}_2(\mathbb{C}))$. We denote this category by $\Rep(U_q(\mathfrak{sl}_2(\mathbb{C})))$. Inside this category is the full monoidal subcategory generated by the $q$ analogue of $\mathbb{C}^2$, which we call $\Fund(U_q(\mathfrak{sl}_2(\mathbb{C})))$. Finally, there is also a $q$ analogue of $\mathcal{TL}$, denoted $\mathcal{TL}_q$ which is the strict pivotal $\mathbb{C}(q)$-linear category generated by a self dual object of dimension $-q - q^{-1}$. Again, we have $\mathcal{TL}_q\cong \Fund(U_q(\mathfrak{sl}_2(\mathbb{C})))$ and as long as $q$ is not a root of unity we also have $\Kar(\Fund(U_q(\mathfrak{sl}_2(\mathbb{C}))))\cong \Rep(U_q(\mathfrak{sl}_2(\mathbb{C})))$.
\end{remark}

\subsection{Main results}
\label{sec-results}

The category $\mathcal{TL}_q$ was generalized by Kuperberg in \cite{Kupe} to describe the representation theory of $U_q(\mathfrak{g})$, when $\mathfrak{g}$ is $\mathfrak{sl}_3, \mathfrak{sp_4}\cong \mathfrak{so}_5$, or $\mathfrak{g_2}$. The only category we consider in this article is the $\mathfrak{sp}_4$ web category, which we denote by $\DD$. We recall the definition here using the convention that $[n]_v:= \dfrac{v^n- v^{-n}}{v-v^{-1}}$ and $[n] := [n]_q\in \mathbb{Z}[q^{\pm 1}]$. 

\begin{defn}
The category $\DD $ is the strict pivotal $\mathbb{Z}[q, q^{-1}, [2]^{-1}]$-linear category generated by two self dual objects, with morphisms generated by 

\begin{equation}
\begin{tikzpicture}[scale=.4, anchorbase]
	\draw[ultra thick, blue] (-1.5,-2) -- (0,0);
	\draw[ultra thick, blue] (1.5, -2) -- (0, 0);
	\draw[ultra thick, green] (0,0) -- (0, 2);
\end{tikzpicture}
\end{equation}

\noindent subject to the following relations.

\begin{equation}\label{circles}
\begin{tikzpicture}[scale=.3, anchorbase]
	\draw[ultra thick, blue] (0,0) circle (1);
\end{tikzpicture}
 \quad = -\dfrac{[6][2]}{[3]} \quad
\end{equation}
\begin{equation}\label{circlet}
\begin{tikzpicture}[scale=.3, anchorbase]
	\draw[ultra thick, green] (0,0) circle (1);
\end{tikzpicture}
 \quad = \dfrac{[6][5]}{[3][2]} \quad
\end{equation}
\begin{equation}\label{monogon}
\begin{tikzpicture}[scale=.3, anchorbase]
	\draw[ultra thick, green] (0, -1) -- (0, 1);
	\draw[ultra thick, blue] (0,2) circle (1);
\end{tikzpicture}
\quad= 0 \quad
\end{equation}
\begin{equation}\label{bigon}
\begin{tikzpicture}[scale=.3, anchorbase]
	\draw[ultra thick, green] (0, -1) -- (0, 1);
	\draw[ultra thick, green] (0, 3) -- (0, 5);
	\draw[ultra thick, blue] (0,2) circle (1);
\end{tikzpicture}
\quad= -[2] \quad
\begin{tikzpicture}[scale=.3, anchorbase]
	\draw[ultra thick, green] (0, -1) -- (0, 5);
\end{tikzpicture}
\end{equation}
\begin{equation}\label{trigon}
\begin{tikzpicture}[scale=.25, anchorbase]
	\draw[ultra thick, green] (-2, 1) -- (-2, 3);
        \draw[ultra thick, green] (2, 1) -- (2, 3);
        \draw[ultra thick, green] (0, -3) -- (0, -5);
        \draw[ultra thick, blue] (-3, 0) -- (-2, 1);
        \draw[ultra thick, blue] (-1, 0) -- (-2, 1);
        \draw[ultra thick, blue] (1, 0) -- (2, 1);
        \draw[ultra thick, blue] (3, 0) -- (2, 1);
        \draw[ultra thick, blue] (-1,0) .. controls (-.55,-1) and (.5,-1) .. (1,0);
        \draw[ultra thick, blue] (-3, 0) -- (-3, -1);
        \draw[ultra thick, blue] (3, 0) -- (3, -1);
        \draw[ultra thick, blue] (-3, -1) -- (0, -3);
        \draw[ultra thick, blue] (3, -1) -- (0, -3);
\end{tikzpicture}
\quad= 0 \quad
\end{equation}
\begin{equation}\label{H=I}
\begin{tikzpicture}[scale=.25, anchorbase]
	\draw[ultra thick, green] (-2, 0) -- (2, 0);
        \draw[ultra thick, blue] (-2, -3) -- (-2, 3);
        \draw[ultra thick, blue] (2, -3) -- (2, 3);
\end{tikzpicture}
= \dfrac{1}{[2]} \quad
\begin{tikzpicture}[scale=.25, anchorbase]
        \draw[ultra thick, blue] (-2, -3) -- (-2, 3);
        \draw[ultra thick, blue] (2, -3) -- (2, 3);
\end{tikzpicture}
\quad+\quad
\begin{tikzpicture}[scale=.25, anchorbase]
	\draw[ultra thick, green] (0, 1) -- (0, -1);
        \draw[ultra thick, blue] (-2, -3) -- (0, -1);
        \draw[ultra thick, blue] (2, -3) -- (0, -1);
        \draw[ultra thick, blue] (-2, 3) -- (0, 1);
        \draw[ultra thick, blue] (2, 3) -- (0, 1);
\end{tikzpicture}
-\dfrac{1}{[2]} \quad
\begin{tikzpicture}[scale=.25, anchorbase]
        \draw[ultra thick, blue] (-2, -3)  .. controls (-1.5,-1) and (1.5,-1) ..  (2, -3);
        \draw[ultra thick, blue] (-2, 3)  .. controls (-1.5,1) and (1.5,1) ..  (2, 3);
\end{tikzpicture}
\end{equation}
\end{defn}

We will refer to these relations as the circle relations, the monogon, bigon, and trigon relation, and the $H\equiv I$ relation. 

Whenever $\ak$ is a field and $q\in \ak$ so that $q+ q^{-1}\ne 0$, we can specialize the diagrammatic category $\DD $ to $\ak$ obtaining $\DDk := \ak \ot \DD$. In \cite[Theorem 5.1, 6.9]{Kupe} Kuperberg proves that the category of finite dimensional (type $1$) representations of $U_q(\mathfrak{sp}_4)$ is equivalent to the Karoubi envelope of $\mathbb{C}(q) \ot \DD$. We write $U_q^{\ak}(\mathfrak{sp}_4)$ to denote the Lusztig divided powers form of the quantum group and $\Tilt(U_q^{\ak}(\mathfrak{sp}_4))$ to denote the monoidal category of tilting modules. It is shown in \cite[Corollary 3.23]{bodish2020web} that $\Tilt(U_q^{\ak}(\mathfrak{sp}_4))$ is equivalent to the Karoubi envelope of $\DDk $. 

Given a tensor product of fundamental representations so that $V(\varpi_1)$ occurs $a$ times and $V(\varpi_2)$ occurs $b$ times, we will write $V(a, b)$ to denote the top summand. In \cite{Kim07}, Kim gives formulas for the $\mathfrak{sp}_4$ clasps projecting to $V(a, 0)$ and $V(0, b)$. The main result of this article is a recursive triple clasp formula for the idempotent projecting to $V(a, b)$. 

\begin{thm}\label{MAINTHM}
Let an oval with $(n, m)$ label denote the idempotent with image $V(a, b)$, then
\begin{equation}
\begin{split}
&\begin{tikzpicture}[scale=.44, anchorbase]
	\draw[ultra thick, black] (-2,-4.2) -- (-2,4.2);
	\draw[ultra thick, blue] (2, -4.2) -- (2, 4.2);
	\draw[ultra thick, black] (1,-4.2) -- (1, 4.2);
	\draw[ultra thick, black, fill= white] (-.3,0) ellipse (2.5 and 1) node{$a+1, b$};
\end{tikzpicture}
\quad = \quad
\begin{tikzpicture}[scale=.44, anchorbase]
	\draw[ultra thick, black] (-2,-4.2) -- (-2,4.2);
	\draw[ultra thick, blue] (2.2, -4.2) -- (2.2, 4.2);
	\draw[ultra thick, black] (1,-4.2) -- (1, 4.2);
	\draw[ultra thick, black, fill= white] (-.55,0) ellipse (2.3 and .8) node{$a, b$};
\end{tikzpicture}
\quad - \dfrac{1}{\kappa_{(a, b), (-1, 1)}}\quad
\begin{tikzpicture}[scale=.44, anchorbase]
	\draw[ultra thick, black] (-2,-4.2) -- (-2,4.2);
	\draw[ultra thick, blue] (2, -4.2) -- (2, -2.4);
	\draw[ultra thick, blue] (2, 2.4) -- (2, 4.2);
	\draw[ultra thick, black] (1,-4.2) -- (1, -3.5);
	\draw[ultra thick, blue] (1,3.5) -- (1, 2.4);
	\draw[ultra thick, black] (1,4.2) -- (1, 3.5);
	\draw[ultra thick, blue] (1,-3.5) -- (1, -2.4);
	\draw[ultra thick, blue] (2,2.4) arc (0:-90:.5);
	\draw[ultra thick, blue] (1,2.4) arc (-180:-90:.5);
	\draw[ultra thick, blue] (2,-2.4) arc (0:90:.5);
	\draw[ultra thick, blue] (1,-2.4) arc (180:90:.5);
	\draw[ultra thick, green] (1.5, -1.9) -- (1.5, 1.9);
	\draw[ultra thick, black, fill= white] (-.55,-3.2) ellipse (2.3 and .8) node{$a, b$};
	\draw[ultra thick, black, fill= white] (-.55,3.2) ellipse (2.3 and .8) node{$a,b$};
	\draw[ultra thick, black, fill= white] (-.3,0) ellipse (2.5 and 1) node{$a-1, b+1$};
\end{tikzpicture}\\
\quad  &- \dfrac{1}{\kappa_{(a, b), (1, -1)}} \quad \begin{tikzpicture}[scale=.44, anchorbase]
	\draw[ultra thick, black] (-2,-4.2) -- (-2, 4.2);
	\draw[ultra thick, blue] (2, -4.2) -- (2, -2.4);
	\draw[ultra thick, blue] (2, 2.4) -- (2, 4.2);
	\draw[ultra thick, black] (1,-4.2) -- (1, -3.5);
	\draw[ultra thick, green] (1,3.5) -- (1, 2.4);
	\draw[ultra thick, black] (1,4.2) -- (1, 3.5);
	\draw[ultra thick, green] (1,-3.5) -- (1, -2.4);
	\draw[ultra thick, blue] (2, 2.4) arc (0:-90:.5);
	\draw[ultra thick, green] (1,2.4) arc (-180:-90:.5);
	\draw[ultra thick, blue] (2,-2.4) arc (0:90:.5);
	\draw[ultra thick, green] (1,-2.4) arc (180:90:.5);
	\draw[ultra thick, blue] (1.5, -1.9) -- (1.5, 1.9);
	\draw[ultra thick, black, fill= white] (-.55,-3.2) ellipse (2.3 and .8) node{$a, b$};
	\draw[ultra thick, black, fill= white] (-.55,3.2) ellipse (2.3 and .8) node{$a,b$};
	\draw[ultra thick, black, fill= white] (-.3,0) ellipse (2.5 and 1) node{$a+1, b-1$};
\end{tikzpicture}
\quad - \dfrac{1}{\kappa_{(a, b), (-1, 0)}}\quad \begin{tikzpicture}[scale=.44, anchorbase]
	\draw[ultra thick, black] (-2,-4.2) -- (-2,4.2);
	\draw[ultra thick, blue] (2, -4.2) -- (2, -2.4);
	\draw[ultra thick, blue] (2, 2.4) -- (2, 4.2);
	\draw[ultra thick, black] (1,-4.2) -- (1, -3.5);
	\draw[ultra thick, blue] (1,3.5) -- (1, 2.4);
	\draw[ultra thick, black] (1,4.2) -- (1, 3.5);
	\draw[ultra thick, blue] (1,-3.5) -- (1, -2.4);
	\draw[ultra thick, blue] (2, 2.4) arc (0:-90:.5);
	\draw[ultra thick, blue] (1,2.4) arc (-180:-90:.5);
	\draw[ultra thick, blue] (2,-2.4) arc (0:90:.5);
	\draw[ultra thick, blue] (1,-2.4) arc (180:90:.5);
	\draw[ultra thick, black, fill= white] (-.55,-3.2) ellipse (2.3 and .8) node{$a, b$};
	\draw[ultra thick, black, fill= white] (-.55,3.2) ellipse (2.3 and .8) node{$a,b$};
	\draw[ultra thick, black, fill= white] (-.3,0) ellipse (2.5 and 1) node{$a-1, b$};
\end{tikzpicture}
\end{split}
\end{equation}
and
\begin{equation}
\begin{split}
&\begin{tikzpicture}[scale=.44, anchorbase]
	\draw[ultra thick, black] (-2,-4.2) -- (-2,4.2);
	\draw[ultra thick, green] (2, -4.2) -- (2, 4.2);
	\draw[ultra thick, black] (1,-4.2) -- (1, 4.2);
	\draw[ultra thick, black, fill= white] (-.3,0) ellipse (2.5 and 1) node{$a, b+1$};
\end{tikzpicture}
\quad = \quad
\begin{tikzpicture}[scale=.44, anchorbase]
	\draw[ultra thick, black] (-2,-4.2) -- (-2,4.2);
	\draw[ultra thick, green] (2.2, -4.2) -- (2.2, 4.2);
	\draw[ultra thick, black] (1,-4.2) -- (1, 4.2);
	\draw[ultra thick, black, fill= white] (-.55,0) ellipse (2.3 and .8) node{$a, b$};
\end{tikzpicture}
\quad - \dfrac{1}{\kappa_{(a, b), (-2, 1)}}\quad
\begin{tikzpicture}[scale=.44, anchorbase]
	\draw[ultra thick, black] (-2,-4.2) -- (-2, 4.2);
	\draw[ultra thick, green] (2, -4.2) -- (2, -2.4);
	\draw[ultra thick, green] (2, 2.4) -- (2, 4.2);
	\draw[ultra thick, black] (1,-4.2) -- (1, -3.5);
	\draw[ultra thick, blue] (1, 3.5) -- (1, 2.4);
	\draw[ultra thick, black] (1,4.2) -- (1, 3.5);
	\draw[ultra thick, blue] (1,-3.5) -- (1, -2.4);
	\draw[ultra thick, green] (2, 2.4) arc (0:-90:.5);
	\draw[ultra thick, blue] (1,2.4) arc (-180:-90:.5);
	\draw[ultra thick, green] (2,-2.4) arc (0:90:.5);
	\draw[ultra thick, blue] (1,-2.4) arc (180:90:.5);
	\draw[ultra thick, blue] (1.5,-1.2) arc (90:180:1);
	\draw[ultra thick, blue] (1.5,1.2) arc (-90:-180:1);
	\draw[ultra thick, black] (.5,-4.2) -- (.5, -3.2);
	\draw[ultra thick, blue] (.5, 3.2) -- (.5, 2.2);
	\draw[ultra thick, black] (.5, 4.2) -- (.5, 3.2);
	\draw[ultra thick, blue] (.5,-3.2) -- (.5, -2.2);
	\draw[ultra thick, blue] (1.5, -1.9) -- (1.5, -1.2);
	\draw[ultra thick, blue] (1.5, 1.9) -- (1.5, 1.2);
	\draw[ultra thick, green] (1.5, -1.2) -- (1.5, 1.2);
	\draw[ultra thick, black, fill= white] (-.55,-3.2) ellipse (2.3 and .8) node{$a, b$};
	\draw[ultra thick, black, fill= white] (-.55,3.2) ellipse (2.3 and .8) node{$a,b$};
	\draw[ultra thick, black, fill= white] (-.3,0) ellipse (2.5 and 1) node{$a-2, b+1$};
\end{tikzpicture}\\
\quad  &- \dfrac{1}{\kappa_{(a, b), (0,0)}} \quad \begin{tikzpicture}[scale=.44, anchorbase]
	\draw[ultra thick, black] (-2,-4.2) -- (-2,4.2);
	\draw[ultra thick, green] (2, -4.2) -- (2, -2.4);
	\draw[ultra thick, green] (2, 2.4) -- (2, 4.2);
	\draw[ultra thick, black] (1,-4.2) -- (1, -3.5);
	\draw[ultra thick, blue] (1,3.5) -- (1, 2.4);
	\draw[ultra thick, black] (1, 4.2) -- (1, 3.5);
	\draw[ultra thick, blue] (1,-3.5) -- (1, -2.4);
	\draw[ultra thick, green] (2, 2.4) arc (0:-90:.5);
	\draw[ultra thick, blue] (1, 2.4) arc (-180:-90:.5);
	\draw[ultra thick, green] (2, -2.4) arc (0:90:.5);
	\draw[ultra thick, blue] (1, -2.4) arc (180:90:.5);
	\draw[ultra thick, blue] (1.5, -1.9) -- (1.5, 1.9);
	\draw[ultra thick, black, fill= white] (-.55,-3.2) ellipse (2.3 and .8) node{$a, b$};
	\draw[ultra thick, black, fill= white] (-.55,3.2) ellipse (2.3 and .8) node{$a,b$};
	\draw[ultra thick, black, fill= white] (-.3,0) ellipse (2.5 and 1) node{$a, b$};
\end{tikzpicture}
\quad - \dfrac{1}{\kappa_{(a, b), (2, -1)}}\quad \begin{tikzpicture}[scale=.44, anchorbase]
	\draw[ultra thick, black] (-2,-4.2) -- (-2,4.2);
	\draw[ultra thick, green] (2, -4.2) -- (2, -1.9);
	\draw[ultra thick, green] (2, 1.9) -- (2, 4.2);
	\draw[ultra thick, black] (1,-4.2) -- (1, -3.5);
	\draw[ultra thick, green] (1,3.5) -- (1, 1.9);
	\draw[ultra thick, black] (1,4.2) -- (1, 3.5);
	\draw[ultra thick, green] (1,-3.5) -- (1, -1.9);
	\draw[ultra thick, blue] (2,-1.9) -- (1, -1.9);
	\draw[ultra thick, blue] (1,1.9) -- (2, 1.9);
	\draw[ultra thick, blue] (1,-1.9) -- (1, 1.9);
	\draw[ultra thick, blue] (2,-1.9) -- (2, 1.9);
	\draw[ultra thick, black, fill= white] (-.55,-3.2) ellipse (2.3 and .8) node{$a, b$};
	\draw[ultra thick, black, fill= white] (-.55,3.2) ellipse (2.3 and .8) node{$a,b$};
	\draw[ultra thick, black, fill= white] (-.3,0) ellipse (2.5 and 1) node{$a+2, b-1$};
\end{tikzpicture}\\
\quad &- \dfrac{1}{\kappa_{(a, b), (0, -1)}} \quad \begin{tikzpicture}[scale=.44, anchorbase]
	\draw[ultra thick, black] (-2,-4.2) -- (-2,4.2);
	\draw[ultra thick, green] (2, -4.2) -- (2, -2.4);
	\draw[ultra thick, green] (2, 2.4) -- (2, 4.2);
	\draw[ultra thick, black] (1,-4.2) -- (1, -3.5);
	\draw[ultra thick, green] (1,3.5) -- (1, 2.4);
	\draw[ultra thick, black] (1,4.2) -- (1, 3.5);
	\draw[ultra thick, green] (1,-3.5) -- (1, -2.4);
	\draw[ultra thick, green] (2,2.4) arc (0:-90:.5);
	\draw[ultra thick, green] (1,2.4) arc (-180:-90:.5);
	\draw[ultra thick, green] (2,-2.4) arc (0:90:.5);
	\draw[ultra thick, green] (1,-2.4) arc (180:90:.5);
	\draw[ultra thick, black, fill= white] (-.55,-3.2) ellipse (2.3 and .8) node{$a, b$};
	\draw[ultra thick, black, fill= white] (-.55,3.2) ellipse (2.3 and .8) node{$a,b$};
	\draw[ultra thick, black, fill= white] (-.3,0) ellipse (2.5 and 1) node{$a, b-1$};
\end{tikzpicture}
\end{split}
\end{equation}
where
\begin{equation}
\kappa_{(a, b), (-1, 1)} = - \dfrac{[a+1]}{[a]}
\end{equation}
\begin{equation}
\kappa_{(a,b), (1, -1)} = \dfrac{[a+2b+3][2b+2]}{[a+ 2b+2][2b]}
\end{equation}
\begin{equation}
\kappa_{(a,b), (-1, 0)} = - \dfrac{[2a+ 2b+ 4][a+ 2b+ 3][a+ 1]}{[2a+ 2b+ 2][a+ 2b+ 2][a]}
\end{equation}
and
\begin{equation}
\kappa_{(a,b), (-2, 1)} = -\dfrac{[a+1][2a+ 2b+4]}{[a-1][2a+ 2b+ 2]}
\end{equation}
\begin{equation}
\kappa_{(a,b), (0,0)} = \dfrac{[a+2][a+2b+4]}{[2][a][a+2b+2]}
\end{equation}
\begin{equation}
\kappa_{(a,b), (2, -1)} = - \dfrac{[2b+2]}{[2b]}
\end{equation}
\begin{equation}
\kappa_{(a,b), (0, -1)} = \dfrac{[2a+2b+4][a+2b+3][2b+2]}{[2a+2b+2][a+2b+1][2b]}
\end{equation}
\end{thm}

Using the double ladder basis for homomorphism spaces in $\DD$ \cite{bodish2020web}, and the ideas from Elias's work on clasps for type $A$ webs \cite{elias2015light}, we can argue that such a recursive formula exists without knowing the $\kappa$'s explicitly. The recursive nature of the clasp formula implies recursive relations among the $\kappa$'s. Our theorem then follows from showing that these relations force the $\kappa$'s to be the values specified in Theorem \eqref{MAINTHM}. 

\subsection{Applications}
\label{sec-applications}

Suppose that $\ak= \mathbb{C}$ and $q= e^{i \pi/\ell}$ for some integer $\ell> 4$. There is a well known construction of a $\mathbb{C}$-linear fusion category as the quotient of $\Tilt(U_q^{\ak}(\mathfrak{sp}_4))$ by the ideal of negligible morphisms \cite{andersen1992}. One may ask for a generators and relations presentation of these quotient categories. Since we have a presentation of the category of tilting modules, it remains to find relations which generate the ideal of negligible morphisms, denoted $\mathcal{N}$. In general the ideal of negligible morphisms is not the monoidal ideal generated by the identity morphisms of all negligible objects. However, the hypothesis $\ell > 4$ guarantees that the non-negligible tilting modules are such that $\dim\Hom(T^{\ak}(\lambda), T^{\ak}(\mu))= \delta_{\lambda, \mu}$. It follows that the quotient by the ideal generated by the negligible objects is a semisimple category. Since the ideal of negligible morphisms is the only monoidal ideal with semisimple quotient, the ideal of negligible morphisms coincides with the monoidal ideal generated by the negligible objects.

There is a $\rho$ shifted and $\ell$ dilated affine Weyl group action on the weight lattice $X$. Extending this action to $\mathbb{R}\ot X$ allows us to partition the dominant weights by their relationship to alcoves in $\mathbb{R}\ot X$. The weights in the interior of the lowest alcove in the cone $- \rho + \mathbb{R}_{\ge 0}\cdot X_+$ are exactly the highest weights of the indecomposable non-negligible tilting modules. Furthermore, every indecomposable negligible tilting module is a direct summand of a tensor product of some tilting module and an indecomposable tilting module with highest weight on the upper closure of the lowest alcove. Thus, the identity morphisms of all indecomposable negligible tilting modules are contained in the monoidal ideal generated by the indecomposable tilting modules on the upper closure of the lowest alcove. 

The tilting modules on the upper closure of the lowest alcove are
\begin{equation}
T^{\ak}\left(2k, \dfrac{\ell-3}{2}-k\right), \ \text{for} \ k= 0, \ldots, \ell-3
\end{equation}
when $\ell$ is odd, and 
\begin{equation}
T^{\ak}\left(k, \dfrac{\ell -4}{2}- k\right), \ \text{for} \ k=0, \ldots, \dfrac{\ell -4}{2}
\end{equation}
when $\ell$ is even. We may conclude that there is a monoidal functor 
\begin{equation}
\DDk \longrightarrow \Tilt(U_q^{\ak}(\mathfrak{sp}_4))/\mathcal{N},
\end{equation}
with kernel the monoidal ideal generated by the $\lambda$ clasps for all $\lambda$ on the upper closure of the lowest alcove. Furthermore, the induced functor
\begin{equation}
\Kar(\DDk /\ker)\longrightarrow \Tilt(U_q^{\ak}(\mathfrak{sp}_4))/\mathcal{N}
\end{equation}
is an equivalence. 

\begin{example}
If $\ell= 5$, then the ideal of negligible morphisms is generated by $\id_{T^{\ak}(0, 1)}$ and $\id_{T^{\ak}(2, 0)}$. The quotient category is equivalent to $\Rep(\mathbb{Z}/2)$, where the sign representation corresponds to $T^{\ak}(1, 0)$.
\end{example}

\begin{example}
If $\ell= 6$, then the ideal of negligible morphisms is generated by $\id_{T^{\ak}(1, 0)}$ and $\id_{T^{\ak}(0, 1)}$. Thus, the quotient category is equivalent to $\bf{Vec}_{\mathbb{C}}$. 
\end{example}

\begin{example}
It is a pleasant exercise to use this approach (along with the diagrammatic category $\mathcal{TL}$) to show that if $\ell= 8$, the negligible quotient of $\Tilt(U_q^{\ak}(\mathfrak{sp}_4))$ is equivalent to $\Tilt(\mathbb{C}(e^{i\pi/4})\ot U_q^{\mathbb{Z}}(\mathfrak{sl}_2))/\mathcal{N}$. Note that if $q= e^{i\pi/8}$, then
\begin{equation}
-\dfrac{[6]_q[2]_q}{[3]_q} = - [2]_{q^2}.
\end{equation}
\end{example}

When $\ell > 4$ is even, the negligible quotient of the category of tilting modules is a modular tensor category \cite[Theorem 4.2]{rowell2006quantumgroupMTC}. One reason this is of interest is that a modular tensor category gives rise to a $321$ dimensional TQFT, and therefore a three manifold invariant \cite[Theorem 3.32]{RTInvariants1991}. In order to compute this invariant explicitly using Kuperberg's graphical calculus for morphisms, it may be useful to have formulas for the non-negligible clasps. The interested reader could consult \cite{KirbyMelvinSL2} to see how this is carried out in the case of the diagrammatic category $\mathcal{TL}$. 

With applications to modular representation theory in mind, there has been some work on writing the idempotents projecting to all indecomposable tilting modules in terms of the double ladder basis for $SL_2(\overline{\mathbb{F}}_p)$ \cite{pJWs} \cite{tubbenhauer2019quivers}. Such idempotents have been referred to as $p$-Jones-Wenzl projectors. Since nobody knows the characters of tilting modules for rank two groups in positive characteristic, this question is not appropriate in our setting. However, the tilting characters are known for the quantum group at a root of unity \cite[Section 8]{SoergelTilt}. A key first step in determining the formulas for the $p$-Jones-Wenzl projectors is to argue that if the characteristic $p$ tilting module is simple, then the characteristic zero clasp can be reduced modulo a maximal ideal to obtain the projector in characteristic $p$. We are careful to point out how this works in the case of $\DD$ \eqref{qJWs}, but do not explore the topic further in this article. 

Lastly, we mention that while preparing the present work for publication, in joint work with Haihan Wu, we solved the analogous problem of finding triple clasp formulas for $\mathfrak{g}_2$ \cite{BodWu}.

\subsection{Structure of the Paper}
\label{sec-struckture}

Section $2$: We recall some facts about the double ladders basis for $\mathfrak{sp}_4$ webs and deduce the triple clasp formula. Section $3$: The recursive formulas for the local intersection forms are stated, and then derived via diagrammatic calculations. We prove the main theorem by showing the conjectured formulas satisfy the recursion. Lastly, we explain how to generalize Elias's clasp conjecture and show that our main theorem verifies this conjecture in type $C_2$.

\subsection{Acknowledgements}
\label{sec-ackackack}

I want to thank Ben Elias for many helpful discussions about idempotents. I also want to thank Haihan Wu and Greg Kuperberg for their suggestion to replace $q^a$ with $A$. The author was supported by NSF grant DMS-1553032.

\section{Clasps}
\label{sec-clasps}

\subsection{Representation theory background}
\label{subsec-repthry}

We briefly recall the notation, conventions, and background found in \cite{bodish2020web}.

Let $\ak$ be a field and let $q\in \ak ^{\times}$ such that $[2]_q \ne 0$. By $\DDk $ we mean the $C_2$ web category base changed to $\ak$. The objects in $\DDk $ are words $\un{w}$ in the alphabet $\lbrace \blues, \greent\rbrace$. 

Let $X= \mathbb{Z}\epsilon_1 \oplus \mathbb{Z}\epsilon_2$ be the set of integral weights for $\mathfrak{sp}_4$. Let $\alpha_1= \epsilon_1- \epsilon_2$ and $\alpha_2= 2\epsilon_2$ be the simple roots for the root system of type $C_2$. Define the dominance order on $X$ such that $\mu < \lambda$ if and only if $\lambda- \mu \in \mathbb{Z}_{\ge 0}\alpha_1 + \mathbb{Z}_{\ge 0}\alpha_2$. Let $X_+$ be the set of dominant integral weights. We write $\varpi_1 = \epsilon_1$ and $\varpi_2 = \epsilon_1+ \epsilon_2$ for the \emph{fundamental weights}, so $X_+ = \lbrace a\varpi_1 + b\varpi_2 \ : \ a, b\in \mathbb{Z}_{\ge 0}\rbrace$. 

Let $U_q^{\mathbb{Z}}(\mathfrak{sp}_4)$ denote Lusztig's divided power form of the quantum group for $\mathfrak{sp}_4$. This algebra has Weyl modules $V^{\mathbb{Z}}(\lambda)$ for all $\lambda\in X_+$. We will write $V^{\ak}(\lambda)$ for $\ak \ot V^{\mathbb{Z}}(\lambda)$. 

An important class of modules are the Weyl filtered modules, which are closed under tensor products. If $X$ is Weyl filtered, we will write $(X: V^{\ak}(\lambda))$ for the filtration multiplicity of $V^{\ak}(\lambda)$ in $X$. This notion is well defined since the classes $[V^{\ak}(\lambda)]$ are a basis for the Grothendieck group of finite dimensional (type $1$) modules. 

In the special case when $\ak = \mathbb{C}$ and $q= 1$ we will write $V(\lambda)$ in place of $V^{\ak}(\lambda)$. The $V(\lambda)$ are the more familiar finite dimensional simple modules of highest weight $\lambda$ for $\mathfrak{sp_4}(\mathbb{C})$. 

There are also dual Weyl modules, which for our purposes can be taken to be $\nabla^{\ak}(\lambda) := V^{\ak}(\lambda)^*$. The modules $\nabla^{\ak}(\lambda)$ and $V^{\ak}(\lambda)$ have the same characters, and
\begin{equation}\label{extvan}
\dim\Ext^i(V(\lambda), \Delta(\mu)) = \delta_{\lambda, \mu}\delta_{i, 0}.
\end{equation}
Moreover, the dual Weyl filtered modules are closed under taking tensor product. Suppose $X$ is a Weyl filtered module and $Y$ is a dual Weyl filtered module, then one can use Equation \eqref{extvan} to show
\begin{equation}\label{dimhomfiltered}
\dim \Hom_{U_q^{\ak}(\mathfrak{sp}_4)}(X, Y)= \sum_{\lambda\in X_+} (X:V^{\ak}(\lambda))(Y: \nabla^{\ak}(\lambda)).
\end{equation}

For any $\lambda\in X_+$, the Weyl modules $V^{\ak}(\lambda)$ fit in between the indecomposable tilting modules $T^{\ak}(\lambda)$, which are filtered by Weyl modules (and by dual Weyl modules) and have $(T^{\ak}(\lambda) : V^{\ak}(\lambda))= 1$, and the simple modules $L^{\ak}(\lambda)$, which are the maximal irreducible quotients of $V^{\ak}(\lambda)$. A Weyl module is irreducible if and only if it is self-dual if and only if it is tilting.

The module $V^{\ak}(\blues) = V^{\ak}(\varpi_1)$ has basis
\[
\lbrace v_{(1,0)}, v_{(-1,1)} = F_sv_{(1, 0)}, v_{(1,-1)} = F_tF_sv_{(1, 0)}, v_{(-1,0)} = F_sF_tF_sv_{(1,0)} \rbrace
\]
and the module $V^{\ak}(\greent) = V^{\ak}(\varpi_2)$ has basis
\[
\lbrace v_{(0,1)} , v_{(2, -1)}= F_tv_{(0,1)}, v_{(0,0)} = F_sF_tv_{(0,1)}, v_{(-2, 1)} = F_s^{(2)}F_tv_{(0,1)}, v_{(0,-1)} = F_tF_s^{(2)}F_tv_{(0,1)}\rbrace.
\]

Our hypothesis that $[2]_q\ne 0$ guarantees that the Weyl modules $V^{\ak}(\varpi_1)$ and $V^{\ak}(\varpi_2)$ are indecomposable tilting modules. The category $\Fund(U_q^{\ak}(\mathfrak{sp}_4))$, is the full subcategory of $\Rep(U_q^{\ak}(\mathfrak{sp}_4))$ monoidally generated by these fundamental tilting modules.  

There is a functor $\eval:\DDk \longrightarrow \Fund(U_q^{\ak}(\mathfrak{sp}_4))$, such that $\eval(\blues) = V^{\ak}(\varpi_1)$ and $\eval(\greent) = V^{\ak}(\varpi_2)$. The functor $\eval$ is is full, faithful, and essentially surjective. 

If $\nu\in X$ and $M$ is a $U_q^{\ak}(\mathfrak{sp}_4)$ module, then we will write $M[\nu]$ for the $\nu$ weight space of $M$. If $M$ is a module which is a direct sum of its weight spaces (note all modules we consider will satisfy this condition), then we will write $\wt M\subset X$ to denote the set of weights of $M$. Given a weight vector $m\in M$, we will write $\wt m$ to denote the weight of the vector. We will also write $\wt \blues = \varpi_1$, $\wt \greent= \varpi_2$, and $\wt \un{w} = \wt w_1+ \ldots + \wt w_n$. 

For $\un{w}= (w_1, \ldots, w_n)$, a word in the alphabet $\lbrace\blues, \greent\rbrace$, we define 
\begin{equation}
V^{\ak}(\un{w}) := V^{\ak}(w_1) \ot \ldots \ot V^{\ak}(w_n),
\end{equation}
and
\begin{equation}
v_{\un{w}, +} := v_{w_1}\ot v_{w_2}\ot \ldots \ot v_{w_n} \in V^{\ak}(\un{w}) 
\end{equation}
where $v_{\blues}= v_{(1, 0)}$ and $v_{\greent}= v_{(0,1)}$. 

Let $P(\un{w})= \lbrace \vec{\nu} = (\nu_1, \ldots \nu_n) \ : \ \nu_i \in \wt V^{\ak}(w_i)\rbrace$. We will write
\begin{equation}
v_{\un{w}, \vec{\nu}} := v_{\nu_1}\ot \ldots \ot v_{\nu_n}\in V^{\ak}(\un{w}). 
\end{equation}
The \emph{subsequence basis} of $V(\un{w})$ is the set
\[
\lbrace v_{\un{w}, \vec{\nu}}\in V^{\ak}(\un{w}) \ : \ \vec{\nu} \in P(\un{w})\rbrace. 
\]

A sequence $(\mu_1, \ldots, \mu_n)$ where $\mu_i\in \wt(V(w_i))$ is a \emph{dominant weight subsequence} of $\un{w}$ if 
\begin{enumerate}
\item $\mu_1$ is dominant. 
\item $V(\mu_1 + \ldots + \mu_{i-1} + \mu_i)$ is a summand of $V(\mu_1 + \ldots + \mu_{i-1})\ot V(w_i)$
\end{enumerate}
If we write $E(\un{w})$ for the set of all dominant weight subsequences of $\un{w}$, then
\begin{equation}
V(\un{w}) \cong \bigoplus_{(\mu_1, \ldots, \mu_n)\in E(\un{w})} V(\mu_1 + \ldots + \mu_n).
\end{equation}
If $\vec{\mu} = (\mu_1, \dots, \mu_n)\in E(\un{w})$, then we will write $\wt \vec{\mu} :=\sum \mu_i$. Therefore, if we denote the multiplicity of $V(\lambda)$ as a summand of $V(\un{w})$ by $[V(\un{w}): V(\lambda)]$ and write
\begin{equation}
E(\un{w}, \lambda):= \lbrace \vec{\mu} \in E(\un{w}) \ : \ \wt \vec{\mu} = \lambda\rbrace,
\end{equation}
then
\begin{equation}
[V(\un{w}): V(\lambda)] = \# E(\un{w}, \lambda).
\end{equation}

\subsection{Recollection of double ladder basis}
\label{subsec-recollectionsaboutdoublelads}

In \cite[Section 2.4, 2.5]{bodish2020web} we defined \emph{elementary light ladders}, \emph{neutral ladders}, \emph{light ladders}, \emph{upside down light ladders}, and \emph{double ladders}. Let us recall the important aspects of these constructions below.

Start by associating an elementary light ladder diagram in $\DD$ to each weight in a fundamental representation: $\mu\mapsto L_{\mu}$. For a dominant weight subsequence $\vec{\mu}\in E(\un{w})$, there are light ladder diagrams $LL_{\un{w}, \vec{\mu}}$. The elementary light ladder diagrams for $\mu_i$ are the building blocks of the light ladder diagram $LL_{\un{w}, \vec{\mu}}$, but in order to make light ladder diagrams out of elementary light ladders we also require neutral diagrams which are used to shuffle words in $\blues$ and $\greent$. The ability to freely choose neutral diagrams means that for a given dominant weight subsequence there may be many choices of light ladder diagrams. 

There is a duality $\mathbb{D}$ on the diagrammatic category, which takes a diagram and flips it upside down. We define upside down light ladders as the image, under $\mathbb{D}$, of usual light ladders.

When defining the double ladder basis, we first fix, for each dominant weight $\lambda$, a choice of object $\un{x}_{\lambda}\in\DD$ satisfying $\wt(\un{x}_{\lambda})= \lambda$. Then for each object $\un{w}$ and each $\vec{\mu}\in E(\un{w})$ we fix a choice of light ladder $LL_{\vec{\mu}}$ such that if $\wt \vec{\mu} = \lambda$, then the target of the diagram is $\un{x}_{\lambda}$. Moreover, for each $\lambda$ there is a unique $\vec{\mu}\in E(\un{x}_{\lambda})$ such that $\wt \vec{\mu} = \lambda$ and we insist the chosen light ladder is the identity. The double ladder basis $\mathbb{LL}$ is then constructed by composing all light ladders with all upside down light ladders.

We refer to the image of these various diagrams under $\eval$ as maps, e.g. $\eval$ applied to a neutral ladder is a neutral map. 

Suppose we have fixed a choice of $\un{x}_{\lambda}$ for all $\lambda \in X_+$ and then fixed a choice of light ladders for all $\un{w}$ and $\vec{\mu}\in E(\un{w})$. We write $\mathbb{LL}$ to denote the associated double ladder basis, $\mathbb{LL}_{\un{w}}^{\un{x}}$ to denote $\mathbb{LL}\cap \Hom_{\DD}(\un{w}, \un{x})$, and $\mathbb{LL}_{\un{w}}^{\un{x}}(\lambda)$ for the collection of double ladders of the form $\mathbb{D}(LL_{\un{x}, \vec{\nu}})\circ LL_{\un{w}, \vec{\mu}}$, where $\vec{\nu} \in E(\un{x}, \lambda)$ and $\vec{\mu}\in E(\un{w}, \lambda)$. Thus, 
\begin{equation}
\mathbb{LL}_{\un{w}}^{\un{x}} = \bigcup_{\lambda\in X_+} \mathbb{LL}_{\un{w}}^{\un{x}}(\lambda). 
\end{equation}
Note that the middle of a diagram in $\mathbb{LL}_{\un{w}}^{\un{x}}(\lambda)$ is the identity of $\un{x}_{\lambda}$. We say that such diagrams \emph{factor through $\lambda$}. 

\begin{defn}
Fix $\lambda\in X_+$. Let $(\DDk)_{< {\lambda}}$ be the $\ak$-linear subcategory whose morphisms are spanned by all double ladders in $\mathbb{LL}_{\un{u}}^{\un{v}}(\chi)$ for all $\un{u}$ and $\un{v}$ and all $\chi< \lambda$.
\end{defn}

\begin{lemma}\label{leftcellular}
Let $f\in \Hom_{\DDk}(\underline{w}, \underline{u})$ and let $\vec{\mu} \in E(\underline{u}, \lambda)$. Then 
\begin{equation}
LL_{\vec{\mu}}\circ f \equiv \sum_{\vec{\nu}\in E(\underline{w}, \lambda)} \ast \cdot LL_{\vec{\nu}} \ \ \ \ \ \text{modulo} \  (\DDk)_{< {\lambda}},
\end{equation}
where $\ast$ represents an element of $\ak$.
\end{lemma}
\begin{proof}
See \cite[Lemma 4.19]{bodish2020web}
\end{proof}

\begin{lemma}\label{cellideal}
Fix $\lambda\in X_+$. The subcategory $(\DDk)_{< {\lambda}}$ is an ideal, i.e. if $D\in (\DDk)_{< {\lambda}}$, then $g\circ D\circ f \in(\DDk)_{< {\lambda}}$ for all $f, g$. 
\end{lemma}
\begin{proof}
See \cite[Lemma 2.8]{ELauda}. Note that $\mathbb{D}$ induces a bijection on homomorphisms spaces in $\DD$ and preserves the basis $\mathbb{LL}$. Thus, the claim follows from Lemma \eqref{leftcellular}.
\end{proof}

\subsection{Downward diagrams and neutral coefficients}
\label{subsec-neutralcoeff}

\begin{defn}
Let $D\in \Hom_{\DDk}(\un{w}, \un{x})$ be an arbitrary diagram. Suppose that there is a horizontal cross section of $D$ which intersects $D$ in the word $\un{y}$. If $\wt \un{y}< \wt \un{w}$, then we say $D$ is a \textbf{downward diagram}. Suppose that there is a horizontal cross section of $D$ which intersects $D$ in the word $\un{y}$. If $\wt \un{y}< \wt \un{x}$, then we say $D$ is an \textbf{upward diagram}. 
\end{defn}

\begin{example}
Any elementary light ladder $L_{\nu}$, for $\nu \notin \lbrace \varpi_1, \varpi_2\rbrace$, is a downward diagram. Any light ladder $LL_{\un{w}, \vec{\mu}}$ is a downward diagram, unless $\mu_i = \wt w_i$ for all $i$. 
\end{example}

\begin{remark}
The duality $\mathbb{D}$ induces a bijection between upward diagrams and downward diagrams. 
\end{remark}

\begin{lemma}\label{strictlydowncell}
If $D\in \Hom_{\DDk}(\un{w}, \un{x})$ is a downward diagram then $D \in (\DDk)_{< \wt \un{w}}$. 
\end{lemma}
\begin{proof}
Suppose that $D$ is a downward diagram, in particular there is some $\un{y}$ with $\wt \un{y} < \wt \un{w}$ such that $D= A\circ B$ for $A\in \Hom_{\DDk}(\un{y}, \un{x})$ and $B\in \Hom_{\DDk}(\un{w}, \un{y})$. If we write $A$ in terms of the double ladder basis, then it is easy to see that $A$ is a linear combination of diagrams in $(\DDk)_{< \wt \un{w}}$. Then from Lemma \eqref{cellideal} it follows that $D\in (\DDk)_{< \wt \un{w}}$. 
\end{proof}

\begin{lemma}\label{downkillstopvector}
Any downward map from $\un{w}$ to $\un{x}$ will send $v_{\un{w}, +}$ to zero.
\end{lemma}
\begin{proof}
Since homomorphisms of $U_q^{\ak}(\mathfrak{sp}_4)$ modules preserve weight spaces, this follows from the observation that if $\wt \un{y} < \wt \un{w}$, then $V^{\ak}(\un{y})[\wt \un{w}] = 0$. 
\end{proof}

\begin{lemma}\label{neutralladders}
If $N: V^{\ak}(\un{w})\rightarrow V^{\ak}(\un{u})$ is a neutral map, then $N(v_{\un{w}, +}) = v_{\un{u}, +}$. Furthermore, if $v_{\un{w}, \vec{\mu}}$ is a subsequence basis element, and $N(v_{\un{w}, \vec{\mu}})$ has a non-zero coefficient for $v_{\un{u}, +}$ after being written in the subsequence basis, then $v_{\un{w}, \vec{\mu}}= v_{\un{w}, +}$. 
\end{lemma}

\begin{proof}
Since $V^{\ak}(\un{w})[\wt \un{w}]=\ak\cdot v_{\un{w}, +}$, the second claim follows from the fact that maps of $U_q^{\ak}(\mathfrak{sp}_4)$ modules preserve weight spaces. Neutral maps are compositions and tensor products of identity maps, $N_{\blues\greent}^{\greent\blues}$, and $N_{\greent\blues}^{\blues\greent}$. So to to verify the first claim we just need to check that $N_{\blues\greent}^{\greent\blues}(v_{\blues\greent, +})= v_{\greent\blues, +}$ and $N_{\greent\blues}^{\blues\greent}(v_{\greent\blues, +})= v_{\blues\greent, +}$. From \cite[Equations 3.27 and 3.30]{bodish2020web} we find
\begin{equation}
N_{\blues\greent}^{\greent\blues}(v_{\blues\greent, +})= L_{(-1, 1)}\ot \id(v_{\blues, +} \ot \mathbb{D}(L_{(-1, 1)})(v_{\greent, +}))= v_{\greent\blues, +}
\end{equation}
and
\begin{equation}
N_{\greent\blues}^{\blues\greent}(v_{\greent\blues, +}) = \id \ot L_{(-1, 1)}(\mathbb{D}(L_{(-1, 1)})(v_{\greent, +})\ot v_{\blues, +})= v_{\blues\greent, +}.
\end{equation}
\end{proof}

\begin{rmk}
In \cite{bodish2020web} we rescaled the generating trivalent vertex in the $B_2$ spider from \cite{Kupe}. Explicitly our trivalent vertex is equal to $\dfrac{1}{\sqrt{[2]}}$ times Kuperberg's trivalent generator. One reason our choice may be preferable to the original, is that using Kuperberg's trivalent vertex the neutral maps have $\xi = [2]$ instead of $\xi = 1$. 
\end{rmk}

Let $\lambda\in X_+$ and let $\un{w}$ and $\un{x}$ be such that $\wt(\un{w})  =\lambda= \wt(\un{x})$. Then $\mathbb{LL}_{\un{w}}^{\un{x}}(\lambda)$ contains a single diagram. Denote this diagram by $I_{\un{w}}^{\un{x}}$. After applying $\eval$, it follows from Lemma \eqref{neutralladders} that $I_{\un{w}}^{\un{x}}$ sends $v_{\un{w}, +}$ to $v_{\un{x}, +}$. 

Suppose we made another choice of $\un{x}_{\lambda}'$ for each dominant weight $\lambda$ (along with choices of all the necessary neutral maps) and then constructed a double ladder basis $\mathbb{LL}'$. Again, there is a unique double ladder diagram $I_{\un{w}}^{'\un{x}}\in \Hom_{\DD}(\un{w}, \un{x})$ which maps $v_{\un{w}, +}$ to $v_{\un{x}, +}$. Since $\mathbb{LL}'$ is a basis, we can express $I_{\un{w}}^{\un{x}}$ as a linear combination of diagrams in $\mathbb{LL}'$
\begin{equation}\label{neutralcoeff}
I_{\un{w}}^{\un{x}} = c\cdot I_{\un{w}}^{'\un{x}} + (\DDk)_{< \wt \un{w}}.
\end{equation}
Looking at how both sides of \eqref{neutralcoeff} act on $v_{\un{w}, +}$ we deduce that $c= 1$. 

\begin{defn}
Let $f\in \Hom_{\DDk}(\un{w}, \un{x})$. We define the \emph{neutral coefficient} of $f$ to be the coefficient of $I_{\un{w}}^{\un{x}}$ when $f$ is expressed in the basis $\mathbb{LL}$. 
\end{defn}

\begin{lemma}
If $\wt \un{w} = \wt \un{x}$, then the neutral coefficient of $f\in \Hom_{\DDk}(\un{w}, \un{x})$ is $c$ if and only if $f(v_{\un{w}, +}) = c\cdot v_{\un{x}, +}$.
\end{lemma}
\begin{proof}
Follows from Lemma \eqref{downkillstopvector} and Lemma \eqref{neutralladders}. 
\end{proof}

\begin{remark}
The discussion given above ensures that the neutral coefficient is independent of any choices that are made in the light ladder algorithm. 
\end{remark}

\subsection{Definition and basic properties of clasps}
\label{subsec-definition}

Our exposition is based on \cite{elias2015light} and \cite[Chapter 11]{soergelbook}.

\begin{defn}
We say that a morphism in $\Hom_{\DDk}(\un{w}, \un{x})$ is a \emph{clasp}, if it is killed by postcomposition with any downward diagram and has neutral coefficient $1$. If $\wt(\un{w}) = \lambda = \wt(\un{x})$, then we may call such a map a $\lambda$-clasp. 
\end{defn}

\begin{lemma}
Let $C\in \Hom_{\DDk}(\un{w}, \un{x})$ have neutral coefficient equal to $1$. Then the following are equivalent:
\begin{enumerate}
\item $C$ is a clasp, 
\item $C$ is killed by postcomposition with any diagram in $(\DDk)_{< \wt \un{w}}$, 
\item $C$ is killed by postcomposition with any diagram of the form $(\id \ot L_{\mu}\ot \id) \circ N$ where $N$ is a neutral diagram and $L_{\mu}$ is an elementary light ladder diagram for $\mu \notin \lbrace \varpi_1, \varpi_2\rbrace$.
\end{enumerate}
\end{lemma}
\begin{proof}
Since a diagram in $(\DDk)_{< \wt \un{w}}$ is a linear combination of downward diagrams, $(1)$ implies $(2)$. The diagram $(\id \ot L_{\mu}\ot \id) \circ N$ in $(3)$ is a downward diagram, thanks to the assumption on $\mu$, so $(1)$ implies $(3)$. From Lemma \eqref{strictlydowncell} we deduce that $(2)$ implies $(1)$. 

By the definition of double ladders as the composition of light ladders and upside down light ladders, and since light ladders are in particular double ladders where the upside down double ladder is the identity, we see that $C$ is killed by postcomposition with $(\DDk)_{< \wt \un{w}}$ if and only if $C$ is killed by postcomposition with any light ladder of the form $LL_{\un{w}, \vec{\mu}}$ where $\wt \vec{\mu}< \wt \un{w}$. Then from the inductive definition of light ladders, we see $(3)$ implies $(2)$. 
\end{proof}

\begin{prop}\label{claspprop}
If a clasp exists then it is unique, and it is also characterized as the map with neutral coefficient $1$ which is killed by precomposition with any upward diagram. The composition of a clasp with a neutral ladder is a clasp (so if any $\lambda$-clasp exists, then all $\lambda$-clasps exist), the composition of two clasps is a clasp, and clasps are preserved by $\mathbb{D}$. 
\end{prop} 
\begin{proof}
We leave it as an exercise to adapt the proof in \cite[Proposition 3.2]{elias2015light} to our setting.
\end{proof}

Graphically we will depict $\lambda$-clasps as ovals labelled by $\lambda$ with source $\un{w}$ and target $\un{x}$. In writing we will denote it by $C_{\lambda}$.

\begin{equation}
C_{\lambda} \quad = \quad \begin{tikzpicture}[scale=.4, anchorbase]
	\draw[ultra thick, black] (-2,-4.2) -- (-2,4.2);
	\draw[ultra thick, black] (2, -4.2) -- (2, 4.2);
	\draw[ultra thick, black, fill= white] (0,0) ellipse (2.5 and 1) node{$\lambda$};
	\node[below] at (0, -4.2) {$\un{w}$};
	\node[above] at (0, 4.2) {$\un{x}$};
\end{tikzpicture}
\end{equation}

The proposition \eqref{claspprop} says that the composition of a clasp with a neutral ladder is a clasp, we will refer to this as \emph{neutral absorption}, depicted diagrammatically as follows. 

\begin{equation}\label{neutralabs}
\begin{tikzpicture}[scale=.4, anchorbase]
	\draw[ultra thick, black] (-2,-4.2) -- (-2,4.2);
	\draw[ultra thick, black] (2, -4.2) -- (2, 4.2);
	\node[below] at (0, -4.2) {$\un{w}$};
	\node[above] at (0, 4.2) {$\un{x}$};
	\draw[ultra thick, blue] (-1, 0) -- (-1, 2.3);
	\draw[ultra thick, green] (-1, 2.3) -- (-1, 4.2);
	\draw[ultra thick, green] (1, 0) -- (1, 2.3);
	\draw[ultra thick, blue] (1, 2.3) -- (1, 4.2);
	\draw[ultra thick, blue] (-1, 2.3) -- (1, 2.3);
	\draw[ultra thick, black, fill= white] (0,0) ellipse (2.5 and 1) node{$\lambda$};
\end{tikzpicture}
\quad = \quad
\begin{tikzpicture}[scale=.4, anchorbase]
	\draw[ultra thick, black] (-2,-4.2) -- (-2,4.2);
	\draw[ultra thick, black] (2, -4.2) -- (2, 4.2);
	\draw[ultra thick, black, fill= white] (0,0) ellipse (2.5 and 1) node{$\lambda$};
	\node[below] at (0, -4.2) {$\un{w}$};
	\node[above] at (0, 4.2) {$\un{x}$};
\end{tikzpicture}
\quad \ \ \ \ \ \ \ \ \ \  \quad
\begin{tikzpicture}[scale=.4, anchorbase]
	\draw[ultra thick, black] (-2,-4.2) -- (-2,4.2);
	\draw[ultra thick, black] (2, -4.2) -- (2, 4.2);
	\node[below] at (0, -4.2) {$\un{w}$};
	\node[above] at (0, 4.2) {$\un{x}$};
	\draw[ultra thick, green] (-1, 0) -- (-1, 2.3);
	\draw[ultra thick, blue] (-1, 2.3) -- (-1, 4.2);
	\draw[ultra thick, blue] (1, 0) -- (1, 2.3);
	\draw[ultra thick, green] (1, 2.3) -- (1, 4.2);
	\draw[ultra thick, blue] (-1, 2.3) -- (1, 2.3);
	\draw[ultra thick, black, fill= white] (0,0) ellipse (2.5 and 1) node{$\lambda$};
\end{tikzpicture}
\quad = \quad
\begin{tikzpicture}[scale=.4, anchorbase]
	\draw[ultra thick, black] (-2,-4.2) -- (-2,4.2);
	\draw[ultra thick, black] (2, -4.2) -- (2, 4.2);
	\draw[ultra thick, black, fill= white] (0,0) ellipse (2.5 and 1) node{$\lambda$};
	\node[below] at (0, -4.2) {$\un{w}$};
	\node[above] at (0, 4.2) {$\un{x}$};
\end{tikzpicture}
\end{equation}

We also observe that the proposition \eqref{claspprop} says the composition of two clasps is a clasp, which is what we will call \emph{clasp absorption}. This is expressed diagrammatically as follows. 

\begin{equation}\label{claspabs}
\begin{tikzpicture}[scale=.4, anchorbase]
	\draw[ultra thick, black] (-2,-4.2) -- (-2,4.2);
	\draw[ultra thick, black] (2, -4.2) -- (2, 4.2);
	\node[below] at (0, -4.2) {$\un{w}$};
	\node[above] at (0, 4.2) {$\un{x}$};
	\draw[ultra thick, black] (-1, 0) -- (-1, 2.3);
	\draw[ultra thick, black] (-1, 2.3) -- (-1, 4.2);
	\draw[ultra thick, black] (1, 0) -- (1, 2.3);
	\draw[ultra thick, black] (1, 2.3) -- (1, 4.2);
	\draw[ultra thick, black, fill= white] (0,0) ellipse (2.5 and 1) node{$\lambda$};
	\draw[ultra thick, black, fill= white] (0,2.3) ellipse (1.5 and .5) node{$\mu$};
\end{tikzpicture}
\quad = \quad 
\begin{tikzpicture}[scale=.4, anchorbase]
	\draw[ultra thick, black] (-2,-4.2) -- (-2,4.2);
	\draw[ultra thick, black] (2, -4.2) -- (2, 4.2);
	\node[below] at (0, -4.2) {$\un{w}$};
	\node[above] at (0, 4.2) {$\un{x}$};
	\draw[ultra thick, black, fill= white] (0,0) ellipse (2.5 and 1) node{$\lambda$};
\end{tikzpicture}
\end{equation}

Finally, note that postcomposing a clasp with a non-identity elementary light ladder results in zero. We will refer to this phenomenon as \emph{clasp orthogonality}, and it can be expressed diagrammatically as follows. 

\begin{equation}\label{clasporth}
\begin{tikzpicture}[scale=.4, anchorbase]
	\draw[ultra thick, black] (-2,-4.2) -- (-2,4.2);
	\draw[ultra thick, black] (2, -4.2) -- (2, 4.2);
	\node[below] at (0, -4.2) {$\un{w}$};
	\node[above] at (0, 4.2) {$\un{x}$};
	\draw[ultra thick, black] (-1, 0) -- (-1, 2.3);
	\draw[ultra thick, black] (-1, 2.3) -- (-1, 4.2);
	\draw[ultra thick, black] (1, 0) -- (1, 2.3);
	\draw[ultra thick, black] (1, 2.3) -- (1, 4.2);
	\draw[ultra thick, black, fill= white] (0,0) ellipse (2.5 and 1) node{$\lambda$};
        \draw[black, ultra thick, fill=white] (-1.5,2) -- (1.5,2) -- (1,3.3) -- (-1,3.3) -- cycle;
        \node at (0, 2.6) {$L_{\nu}$};
        \end{tikzpicture}
\quad = 0 \quad 
\end{equation}

\begin{remark}
The $\lambda$ clasps give a \emph{compatible system of idempotents} \cite[Definition 3.3]{elias2015light}, and therefore represent an object in the Karoubi envelope of $\DDk $. This object is a common summand of the objects $\un{w}$ such that $\wt(\un{w})= \lambda$. 
\end{remark}

\begin{notation}
Given an idempotent $e\in \End_{\ak\ot \DD}(\un{w})$ we get an object $(\un{w}, e)$ in $\text{Kar}(\DDk )$. The object in the Karoubi envelope which corresponds to $(\un{w}, C_{\lambda})$ will be denoted by $\lambda$, for all $\un{w}$ such that $\wt \un{w} = \lambda$. 
\end{notation}

Recall that in the Karoubi envelope we have
\begin{equation}
\Hom_{\text{Kar}(\DDk )}((\un{w}, e), (\un{x}, f))= f\Hom_{\ak\ot \DD}(\un{w}, \un{x})e. 
\end{equation}

\begin{cor}\label{scalarend}
Suppose the $\lambda$ and $\chi$ clasps both exist. Then $\Hom_{\text{Kar}(\DDk )}(\lambda, \chi)$ is spanned by the identity if $\lambda = \chi$ and is zero otherwise. 
\end{cor}
\begin{proof}
Compare with \cite[Corollary 3.6]{elias2015light}. Let $D\in \mathbb{LL}_{\un{x}_{\lambda}}^{\un{x}_{\chi}}$. From clasp orthogonality it follows that $C_{\chi} \circ D \circ C_{\lambda}= 0$ unless $D= I_{\un{x}_{\lambda}}^{\un{x}_{\chi}}$. If $\lambda \ne \chi$, then every diagram in $\mathbb{LL}_{\un{x}_{\lambda}}^{\un{x}_{\un{\chi}}}$ is strictly lower, so $\Hom_{\text{Kar}(\DDk )}(\lambda, \chi)= 0$.  Thanks to neutral absorption \eqref{neutralabs} we have $C_{\chi}\circ I_{\un{x}_{\lambda}}^{\un{x}_{\chi}} \circ C_{\lambda}= C_{\chi} \circ \id \circ C_{\lambda}$. Then from clasp absorption \eqref{claspabs} it follows that $C_{\chi} = C_{\chi} \circ \id \circ C_{\lambda} = C_{\lambda}$.  
\end{proof}

\begin{lemma}
The $\lambda$ clasp exists in $\DDk $ if and only if $V^{\ak}(\lambda)$ is a direct summand of $V^{\ak}(\un{x}_{\lambda})$. Moreover, when the $\lambda$ clasp exists we have $C_{\lambda} = \eval^{-1}(e_{\lambda})$ where $e_{\lambda}$ is the endomorphism of $V^{\ak}(\un{x}_{\lambda})$ projecting to $V^{\ak}(\lambda)$. 
\end{lemma}

\begin{proof}
Suppose that the $\lambda$ clasp does exist. Consider the idempotent $e_{\lambda}\in \End(V^{\ak}(\un{x}_{\lambda}))$ which is the image under $\eval$ of the $\lambda$ clasp in $\End_{\DDk}(\un{x}_\lambda)$. The map $e_{\lambda}$ projects to a direct summand of $V^{\ak}(\un{x}_{\lambda})$, and by Corollary \eqref{scalarend} the summand has endomorphism ring $\ak \cdot \id$. Since $V^{\ak}(\un{x}_{\lambda})[\lambda] = \ak \cdot v_{\un{x}_{\lambda}, +}$ and the lambda clasp preserves the $\lambda$ weight vector $v_{\un{x}_{\lambda}, +}$, the object $\im(e_{\lambda})$ has a one dimensional lambda weight space. An object with a one-dimensional $\lambda$ weight space and a local endomorphism ring must be the indecomposable tilting module of highest weight $\lambda$. Since the endomorphism ring of a tilting module is $\ak \cdot \id$ if and only if the indecomposable tilting module is an irreducible Weyl module, it follows that $\im(e_{\lambda})\cong V^{\ak}(\lambda)$. 

Suppose $V^{\ak}(\lambda)$ is a summand of $V^{\ak}(\un{x}_{\lambda})$, so there is an idempotent $e_{\lambda}\in \End_{U_q^{\ak}(\mathfrak{sp}_4)}(V^{\ak}(\un{x}_{\lambda}))$ which projects to $V^{\ak}(\lambda)$. Since $V^{\ak}(\lambda)[\lambda]$ is one dimensional, it follows that $\im(e_{\lambda})[\lambda] = V^{\ak}(\un{x}_{\lambda})[\lambda] = \ak \cdot v_{\un{x}_{\lambda, +}}$. Restricting $e_{\lambda}$ to the $\lambda$ weight space induces an isomorphism. Hence, $e_{\lambda}(v_{\un{x}_{\lambda}, +})= \xi v_{\un{x}_{\lambda}, +}$ for some $\xi \in \ak^{\times}$. Since $e_{\lambda}$ is idempotent, $\xi =1$, so $e_{\lambda}$ has neutral coefficient one. By Lemma \eqref{downkillstopvector}, postcomposing $e_{\lambda}$ with a downward map induces a map $V^{\ak}(\lambda) \rightarrow V^{\ak}(\un{y})$ which has $V^{\ak}(\lambda)[\lambda]$ in its kernel, and therefore is zero. We conclude that $e_{\lambda}$ is a clasp. 
\end{proof}

\begin{remark}
Since the finite dimensional representations of $\mathfrak{sp}_4(\mathbb{C})$ are completely reducible, it follows that when $\ak = \mathbb{C}$ and $q= 1$, $\lambda$ clasps exist for all $\lambda \in X_+$. In more generality, if $K$ is any field and $q$ is transcendental, then $\lambda$ clasps exist over $K$ for all $\lambda\in X_+$. 
\end{remark}

\begin{remark}\label{qJWs}
We argue that if $V^{\ak}(\lambda)$ is a direct summand of $V^{\ak}(\un{x}_{\lambda})$, i.e. the $\lambda$ clasp exists over $\ak$, then the characteristic zero clasp can be used to compute the $\ak$ clasp. 

Let $\mathcal{A}= \mathbb{Z}[q, q^{-1}, [2]_q^{-1}]$. When we say ``all fields" we mean all pairs $\ak$ and $q\in \ak$ such that $q+ q^{-1}\ne 0$. Any quotient of $\mathcal{A}$ by a maximal ideal will give such a pair. 

From \cite{Kupe} we know that the set $\textbf{D}$ of non-elliptic webs spans $\DD$ over $\mathcal{A}$, and we know from \cite{bodish2020web} that $\mathbb{LL}$ is linearly independent over all fields $\ak$. It follows that the coefficients of a linear dependence among double ladders over $\mathcal{A}$ must all be contained in every maximal ideal of $\mathcal{A}$. But the Jacobson radical of $\mathcal{A}$, denoted $J(\mathcal{A})$, is zero. So $\mathbb{LL}$ is linearly independent over $\mathcal{A}$. Furthermore, we established in \cite{bodish2020web} that the sets $\textbf{D}$ and $\mathbb{LL}$ both give bases of $\DDk $ for all fields $\ak$. 

Fix objects $\un{w},\un{u}\in \DD$. Consider the matrix which expresses a double ladder in terms of the spanning set $\textbf{D}$ 
\begin{equation}
A_{\un{w}}^{\un{u}}: \mathcal{A}\mathbb{LL}_{\un{w}}^{\un{u}}\longrightarrow \mathcal{A}B_{\un{w}}^{\un{u}} =\Hom_{\DD}(\un{w}, \un{u}). 
\end{equation}
This matrix is an isomorphism over $\ak$, for all fields $\ak$, so the determinant is not contained in any maximal ideal in $\mathcal{A}$. Thus, $\det A_{\un{w}}^{\un{u}}$ is a unit in $\mathcal{A}$, and $A_{\un{w}}^{\un{u}}$ is invertible over $\mathcal{A}$. Hence, $\mathbb{LL}$ spans $\DD$ over $\mathcal{A}$. 

Let $\mathcal{O}$ be a complete discrete valuation ring, which is an $\mathcal{A}$ algebra, and such that $\mathcal{O}/\mathfrak{m} = \ak$. Assume that the field $K = \text{Frac}(\mathcal{O})$ is characteristic zero and $q\in K$ is transcendental. Suppose that $V^{\ak}(\lambda)$ is a summand of $V^{\ak}(\un{x}_{\lambda})$, i.e. there is a clasp $e_{\lambda}^{\ak}\in \End_{\ak\ot \DD}(\un{x}_{\lambda})$. The endomorphism $e_{\lambda}^{\ak}$ is an idempotent so it can be lifted to $e_{\lambda}^{\mathcal{O}}\in \End_{\mathcal{O}\ot \DD}(\un{x}_{\lambda})$. Since $\mathbb{LL}$ is a basis of $\DD$ over $\mathcal{A}$, it follows that $\mathbb{LL}$ is a basis of the category $\mathcal{O}\ot \DD$. Therefore, 
\[
e_{\lambda}^{\mathcal{O}}= x\id + (\mathcal{O}\ot \DD)_{< \lambda}.
\] 
Since $e_{\lambda}^{\mathcal{O}}$ specializes to $e_{\lambda}^{\ak}$, which in turn sends $v_{\un{x}_{\lambda}, +}$ to $v_{\un{x}_{\lambda}, +}$, there is some $m\in \mathfrak{m}$ such that $x= 1-m$. Also, $e_{\lambda}^{\mathcal{O}}$ is an idempotent and $(\mathcal{O}\ot \DD)_{< \lambda}$ is an ideal in $\mathcal{O}\ot \DD$ so
\[
x\id + (\mathcal{O}\ot \DD)_{< \lambda} = \left(x\id + (\mathcal{O}\ot \DD)_{< \lambda}\right)^2 = x^2\id + (\mathcal{O}\ot \DD)_{< \lambda}.
\]
Comparing neutral coefficients, we find $x= x^2$. It follows that $1- m= (1-m)^2=1- 2m + m^2$, which implies $m^2 = m$. Since $m\in \mathfrak{m}$, we may conclude that $m= 0$. 

From the fact that $\mathbb{LL}$ is a basis over $\mathcal{O}$ it follows that the homomorphism spaces in $\mathcal{O}\ot \DD$ are free and finitely generated $\mathcal{O}$-modules. Thus, the $\mathcal{O}$ module 
\begin{equation}
e^{\mathcal{O}}_{\lambda}\Hom_{\mathcal{O}\ot \DD}(\un{x}_{\lambda})e^{\mathcal{O}}_{\lambda}
\end{equation}
is a finitely generated projective $\mathcal{O}$ module. Since $\mathcal{O}$ is local, one can use Nakayama's lemma to show that projective and finitely generated implies free of finite rank. A consequence is the equality
\begin{equation}
\text{rk}_{\mathcal{O}} \ e^{\mathcal{O}}_{\lambda}\End_{\mathcal{O}\ot \DD}(\un{x}_{\lambda})e^{\mathcal{O}}_{\lambda}= \dim_{\ak}e^{\ak}_{\lambda}\End_{\ak\ot \DD}(\un{x}_{\lambda})e^{\ak}_{\lambda}. 
\end{equation}
We know $e^{\ak}_{\lambda}\End_{\ak\ot \DD}(\un{x}_{\lambda})e^{\ak}_{\lambda} = \ak \cdot e^{\ak}_{\lambda}$, so we may deduce that $e^{\mathcal{O}}_{\lambda}\Hom_{\mathcal{O}\ot \DD}(\un{x}_{\lambda})e^{\mathcal{O}}_{\lambda} = \mathcal{O}\cdot e^{\mathcal{O}}_{\lambda}$.

On the other hand, we know there is a characteristic zero clasp, $e_{\lambda}^K\in \End_{K\ot \DD}(\un{x}_{\lambda})$. Using that $e_{\lambda}^K$ has neutral coefficient one and is orthogonal to all downward diagrams, we may conclude that $e_{\lambda}^{\mathcal{O}}e_{\lambda}^Ke_{\lambda}^{\mathcal{O}}= e_{\lambda}^K$, so $e_{\lambda}^K \in e_{\lambda}^{\mathcal{O}}\End_{K\ot \DD}(\un{x}_{\lambda})e_{\lambda}^{\mathcal{O}}= Ke_{\lambda}^{\mathcal{O}}$. By comparing neutral coefficients we see that $e_{\lambda}^K = e_{\lambda}^{\mathcal{O}}$. 

Over $K$ the $\lambda$ clasp exists for all $\lambda\in X_+$, so to compute $e^K_{\lambda}$ we are free to use the recursion from our main theorem to expand this clasp in terms of the basis $\mathbb{LL}$. The argument we just sketched implies that the coefficients, of the double ladders, in the expanded clasp actually lie in $\mathcal{O}$. So we can reduce $e^{K}_{\lambda}$ modulo a maximal ideal to obtain $e^{\ak}_{\lambda}$. 
\end{remark}



\subsection{Intersection forms and triple clasp formulas}
\label{subsec-tripleclasp}

Let $a\in \lbrace \blues, \greent\rbrace$. If $\ak= \mathbb{C}$ and $q= 1$, we know that $\eval(\lambda\ot a) = V(\lambda)\ot V(\varpi_a)$ decomposes as described by $(2.6a)$ and $(2.6b)$ in \cite{bodish2020web}. 

\begin{defn}
Let $a\in \lbrace \blues, \greent\rbrace$ and let $\lambda \in X_+$. Define the set $S_{\lambda, a}$ to be the collection of weights $\mu \in \wt V(a)$ such that $V(\lambda+ \mu)$ is a direct summand of $V(\lambda) \ot V(a)$.
Since each weight in $\wt V(a)$ is multiplicity one,
\begin{equation}\label{charzerodecomp}
V(\lambda) \ot V(a) \cong \bigoplus_{\mu\in S_{\lambda, a}} V(\lambda+ \mu).
\end{equation}
\end{defn}

\begin{lemma}
$V^{\ak}(\lambda)\ot V^{\ak}(a)$ has a filtration by the Weyl modules $V^{\ak}(\lambda+ \mu)$ for $\mu \in S_{\lambda,a}$.
\end{lemma}
\begin{proof}
The tensor product of Weyl filtered modules has a Weyl module filtration. Since $V^{\ak}(\chi)$ has the same character as $V(\chi)$, the filtration multiplicities are determined by the character of the Weyl filtered module. Therefore, the claim follows from Equation \eqref{charzerodecomp}. 
\end{proof}

\begin{defn}
Let $\mu\in S_{\lambda, a}$. Suppose the $\lambda$ and $\lambda+ \mu$ clasps exist over $\ak$. There is an elementary light ladder $L_{\mu}$ which induces a map 
\[
C_{\lambda + \mu} \circ (\id \ot L_{\mu}) \circ (C_{\lambda} \ot \id_a) : \lambda \ot a\to \lambda+ \mu.
\] 
We denote this map by $E_{\lambda, \mu}$ and depict it diagrammatically by

\begin{equation}\label{Emap}
E_{\lambda, \mu}= \quad
\begin{tikzpicture}[scale=.4, anchorbase]
	\draw[ultra thick, black] (-3,-4.2) -- (-3,4.2);
	\draw[ultra thick, black] (3, -4.2) -- (3, 4.2);
	\draw[ultra thick, black] (0, -4.2) -- (0, 4.2);
	\node[below] at (3, -4.2) {$a$};
	\draw[ultra thick, black, fill= white] (-1.5,-2.8) ellipse (3 and 1) node{$\lambda$};
	\draw[ultra thick, black, fill= white] (0,2.8) ellipse (4 and 1) node{$\lambda + \mu$};
	\draw[black, ultra thick, fill=white] (-1.5,-1) -- (4.5,-1) -- (3.5,1) -- (-.5, 1) -- cycle;
        \node at (1.5, 0) {$L_{\mu}$};
\end{tikzpicture}
\end{equation}
\end{defn}

\begin{prop}\label{claspedbasis}
Suppose the $\lambda$ clasp exists and that the $\lambda +\mu$ clasp exists. Then $\lbrace E_{\lambda, \mu}\rbrace$ is a basis for $\Hom_{\text{Kar}\DDk }(\lambda\ot a, \lambda+\mu)$.
\end{prop}

\begin{proof} 
Since double ladders are a basis, it follows that after postcomposing with $C_{\lambda+ \mu}$ and precomposing with $C_{\lambda}\ot \id_a$ the double ladders $ \mathbb{LL}_{\un{x}_{\lambda} a}^{\un{x}_{\lambda + \mu}}$ will span $\Hom_{\text{Kar}\DDk }(\lambda\ot a, \mu)$.

Let $D\in \mathbb{LL}_{\un{x}_{\lambda} a}^{\un{x}_{\lambda + \mu}}$. By the definition of double ladders, there are dominant weight sequences $\vec{\nu}\in E(\un{x}_{\lambda+ \mu})$ and $\vec{\chi}= (\chi_1, \ldots \chi_n)\in E(\un{x}_{\lambda}a)$ such that $D= \mathbb{D}(LL_{\un{x}_{\lambda+ \mu}, \vec{\nu}})\circ LL_{\un{x}_{\lambda}a, \vec{\chi}}$. Due to clasp orthogonality \eqref{clasporth}, $C_{\lambda+ \mu} \circ D \circ (C_{\lambda}\ot \id_a) = 0$ unless $\vec{\nu}\in E(\un{x}_{\lambda+ \mu}, \lambda+ \mu)$ and $(\chi_1, \ldots \chi_{n-1})\in E(\un{x}_{\lambda}, \lambda)$. Using the neutral absorption property of clasps \eqref{neutralabs}, we now see that 
\[
E_{\lambda, \mu} = C_{\lambda + \mu} \circ \mathbb{LL}_{\un{x}_{\lambda} a}^{\un{x}_{\lambda+ \mu}} \circ (C_{\lambda} \ot \id_a).
\]
Therefore, $E_{\lambda, \mu}$ spans $\Hom_{\text{Kar}\DDk }(\lambda\ot a, \lambda + \mu)$. 

The $\lambda+ \mu$ clasp exists, so $V^{\ak}(\lambda + \mu)$ is an irreducible Weyl module, and therefore $V^{\ak}(\lambda+ \mu)\cong V^{\ak}(\lambda+ \mu)^*$. Since $V^{\ak}(\lambda+\mu)$ occurs exactly once in the Weyl filtration of $V^{\ak}(\lambda) \ot V^{\ak}(a)$, it follows from Equation \eqref{dimhomfiltered} that
\[
\dim \Hom_{\text{Kar}\DDk }(\lambda\ot a, \lambda + \mu) = \dim \Hom_{U_q^{\ak}(\mathfrak{sp}_4)}(V^{\ak}(\lambda)\ot V^{\ak}(a), V^{\ak}(\lambda+ \mu)^*) = 1.
\]
Thus $E_{\lambda, \mu}$ spans a one dimensional vector space and therefore is a basis. 
\end{proof}


\begin{defn}
The map $K_{\lambda, \mu} := E_{\lambda, \mu}\circ \DM E_{\lambda, \mu}$ is an endomorphism of $\lambda + \mu$, and this endomorphism space is spanned by the identity map. We define the local intersection form $\kappa_{\lambda, \mu}$ to be the neutral coefficient of $E_{\lambda, \mu}\circ \DM E_{\lambda, \mu}$. 

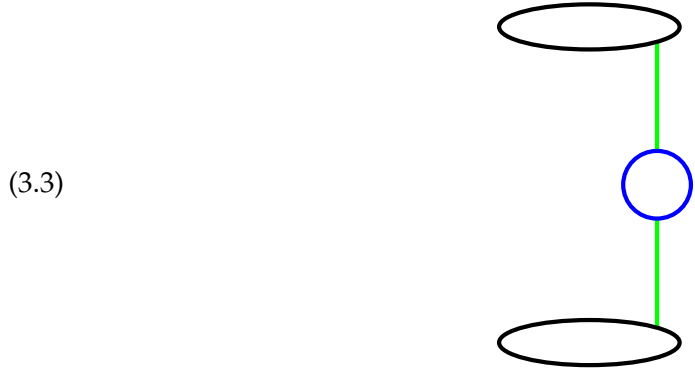
\begin{equation}\label{LIFdefn1}
K_{\lambda, \mu}= \quad
\begin{tikzpicture}[scale=.3, anchorbase]
	\draw[ultra thick, black] (-3, 0) -- (-3,8.4);
	\draw[ultra thick, black] (3, 0) -- (3, 8.4);
	\draw[ultra thick, black] (0, 0) -- (0, 8.4);
	\draw[ultra thick, black, fill= white] (-1.5,1.4) ellipse (3 and 1) node{$\lambda$};
	\draw[ultra thick, black, fill= white] (0,7) ellipse (4 and 1) node{$\lambda + \mu$};
	\draw[black, ultra thick, fill=white] (-1.5,3.2) -- (4.5,3.2) -- (3.5,5.2) -- (-.5, 5.2) -- cycle;
        \node at (1.5, 4.2) {$L_{\mu}$};
        \draw[ultra thick, black] (-3,0) -- (-3, -8.4);
	\draw[ultra thick, black] (3, 0) -- (3, -8.4);
	\draw[ultra thick, black] (0, 0) -- (0, -8.4);
	\draw[ultra thick, black, fill= white] (-1.5,-1.4) ellipse (3 and 1) node{$\lambda$};
	\draw[ultra thick, black, fill= white] (0,-7) ellipse (4 and 1) node{$\lambda + \mu$};
	\draw[black, ultra thick, fill=white] (-1.5,-3.2) -- (4.5,-3.2) -- (3.5,-5.2) -- (-.5, -5.2) -- cycle;
        \node at (1.5, -4.2) {$\mathbb{D}(L_{\mu})$};
\end{tikzpicture}
\quad = \kappa_{\lambda, \mu} \quad
\begin{tikzpicture}[scale=.3, anchorbase]
	\draw[ultra thick, black] (-3, -8.4) -- (-3,8.4);
	\draw[ultra thick, black] (3, -8.4) -- (3, 8.4);
	\draw[ultra thick, black] (0, -8.4) -- (0, 8.4);
	\draw[ultra thick, black, fill= white] (0,0) ellipse (4 and 1) node{$\lambda + \mu$};
\end{tikzpicture}
\end{equation}
\end{defn}

\begin{lemma}\label{whenissplit}
Suppose that both the $\lambda$ clasp and the $\lambda+ \mu$ clasp exist. If the local intersection form $\kappa_{\lambda, \mu}$ is nonzero, then $\dfrac{1}{\kappa_{\lambda, \mu}}\mathbb{D}E_{\lambda, \mu} \circ E_{\lambda, \mu}$ is an idempotent in $\End_{\text{Kar}(\DDk )}(\lambda \ot a)$ which projects to $V^{\ak}(\lambda+ \mu)$. If the local intersection form is zero, then $V^{\ak}(\lambda+ \mu)$ is not a summand of $V^{\ak}(\lambda)\ot V^{\ak}(a)$. 
\end{lemma}

\begin{proof}
If $\kappa_{\lambda, \mu} \ne 0$, then $\dfrac{1}{\kappa_{\lambda, \mu}}\mathbb{D}E_{\lambda, \mu} \circ E_{\lambda, \mu}$ is a non-zero idempotent factoring through $V^{\ak}(\lambda+ \mu)$. Since the $\lambda + \mu$ clasp exists, the module $V^{\ak}(\lambda+ \mu)$ is irreducible. Thus, the idempotent has image isomorphic to $V^{\ak}(\lambda+ \mu)$. 

Since $E_{\lambda, \mu}$ is a basis for $\Hom_{\DDk}(\lambda\ot a, \lambda+ \mu)$ it follows that $\mathbb{D}(E_{\lambda, \mu})$ is a basis for $\Hom_{\DDk}(\lambda+ \mu, \lambda \ot a)$. If $\kappa_{\lambda, \mu} = 0$, then every pair of projection $V^{\ak}(\lambda)\ot V^{\ak}(a)\rightarrow V^{\ak}(\lambda+ \mu)$ and inclusion $V^{\ak}(\lambda+ \mu) \rightarrow V^{\ak}(\lambda)\ot V^{\ak}(a)$ compose to be $0 \in \End_{U_q^{\ak}(\mathfrak{sp}_4)}(V^{\ak}(\lambda+ \mu))$. Thus, $V^{\ak}(\lambda+ \mu)$ is not a direct summand of $V^{\ak}(\lambda) \ot V^{\ak}(a)$. 
\end{proof}

\begin{remark}
By working with the ideals $(\DDk)_{< \lambda}$ instead of with clasps, one can show that the indecomposable tilting module $T^{\ak}(\lambda + \mu)$ is a direct summand of $V^{\ak}(\lambda) \ot V^{\ak}(a)$ if and only if $\kappa_{\lambda, \mu} \ne 0$. 
\end{remark}

\begin{prop}\label{claspformula}
Suppose that the $\lambda$ clasp exists. Also assume the $\lambda + \mu$ clasps exist and the $\kappa_{\lambda, \mu}$ are invertible in $\ak$, for all $\mu\in S_{\lambda, a} - \lbrace \varpi_a\rbrace$. Then the $\lambda + \varpi_a$ clasp exists and 
\begin{equation}\label{clasprecurs}
C_{\lambda+ \varpi_a} = C_{\lambda} \ot \id_a - \sum \kappa_{\lambda, \mu}^{-1} \mathbb{D}E_{\lambda, \mu} \circ E_{\lambda, \mu},
\end{equation}
where the sum is over all $\mu\in S_{\lambda, a} - \lbrace \varpi_a\rbrace$. 
\end{prop}

\begin{proof}
Lemma \eqref{whenissplit} implies that $V^{\ak}(\lambda+ \mu)$ is a direct summand of $V^{\ak}(\lambda) \ot V^{\ak}(a)$ for all $\mu\in S_{\lambda, a}- \lbrace \varpi_a \rbrace$. So
\begin{equation}\label{decomptensor}
V^{\ak}(\lambda) \ot V^{\ak}(a) = X \bigoplus_{\mu \in S_{\lambda, a}- \lbrace \varpi_a \rbrace} V^{\ak}(\lambda + \mu), 
\end{equation}
where $X$ is a direct summand with highest weight $\lambda + \varpi_a$. Direct summands of Weyl filtered modules have a Weyl filtration, so $X$ is filtered by Weyl modules. Comparing characters of both sides of \eqref{decomptensor} implies that $X \cong V^{\ak}(\lambda+ \varpi_a)$. 
The decomposition in \eqref{decomptensor} gives rise to the following equality in the endomorphism algebra of $V^{\ak}(\lambda)\ot V^{\ak}(a)$:
\begin{equation}\label{endoequality}
\id_{V^{\ak}(\lambda)} \ot \id_{V^{\ak}(a)} = \sum_{\mu \in S_{\lambda, a}} e_{\lambda + \mu},
\end{equation}
where $e_{\lambda + \mu}$ is the projection to the summand isomorphic to $V^{\ak}(\lambda + \mu)$. The module $V^{\ak}(\lambda)$ is a direct summand of $V^{\ak}(\un{x}_{\lambda})$, the image of $C_{\lambda}$, so $V^{\ak}(\lambda)\ot V^{\ak}(a)$ is a direct summand of $V^{\ak}(\un{x}_{\lambda}a)$, the image of $C_{\lambda}\ot \id_a$. Precomposing and postcomposing Equation \eqref{endoequality} with $C_{\lambda}\ot \id_a$ yields the desired equality from Equation \eqref{clasprecurs}, in the endomorphism algebra of $V^{\ak}(\un{x}_{\lambda}a)$.
\end{proof}

\section{Clasp coefficients}
\label{subsec-coeffs}
\subsection{Deriving the recursive formulas}
\label{subsec-deriverecursive}

We will compute recursive formulas for the local intersection forms $\kappa_{\lambda, \mu}$ using the graphical calculus for $\DD$. 

\begin{notation}
To simplify notation, we will write $(a,b)$ for $a\varpi_1 + b\varpi_2$. We will also leave off labels of clasps when the highest weight is understood. Furthermore, we will often leave off extra strands below (above) clasps which are on the bottom (top) of the diagram, as well as strands to the left of a diagram which has a clasp at the top or bottom. This is justified because all clasps with the same highest weight are transformed to one another by applying neutral diagrams on the top and bottom (in other contexts this could be nontrivial to verify, but it is easy to see that any two words in $\blues$ and $\greent$ of the same weight differ by a neutral diagram). We also freely use clasp absorption \eqref{claspabs} to simplify formulas. For example \eqref{LIFdefn1} becomes

\begin{equation}\label{LIFdefn2}
K_{\lambda, \mu}= \quad
\begin{tikzpicture}[scale=.3, anchorbase]
	\draw[ultra thick, black] (3, 0) -- (3, 7.5);
	\draw[ultra thick, black] (0, 0) -- (0, 7.5);
	\draw[ultra thick, black, fill= white] (0,7) ellipse (4 and 1);
	\draw[black, ultra thick, fill=white] (-1.5,3.2) -- (4.5,3.2) -- (3.5,5.2) -- (-.5, 5.2) -- cycle;
        \node at (1.5, 4.2) {$L_{\mu}$};
	\draw[ultra thick, black] (3, 0) -- (3, -7.5);
	\draw[ultra thick, black] (0, 0) -- (0, -7.5);
	\draw[ultra thick, black, fill= white] (-1.5,0) ellipse (3 and 1);
	\draw[ultra thick, black, fill= white] (0,-7) ellipse (4 and 1);
	\draw[black, ultra thick, fill=white] (-1.5,-3.2) -- (4.5,-3.2) -- (3.5,-5.2) -- (-.5, -5.2) -- cycle;
        \node at (1.5, -4.2) {$\mathbb{D}(L_{\mu})$};
\end{tikzpicture}
\quad = \kappa_{\lambda, \mu} \quad
\begin{tikzpicture}[scale=.3, anchorbase]
	\draw[ultra thick, black, fill= white] (0,0) ellipse (4 and 1);
\end{tikzpicture}
\end{equation}
and
\begin{equation}
\begin{tikzpicture}[scale=.3, anchorbase]
        \draw[ultra thick, black] (-3, -8.2) -- (-3, 8.2);
        \draw[ultra thick, black] (-1.5, -8.2) -- (-1.5, 8.2);
        \draw[ultra thick, black] (0, -8.2) -- (0, 8.2);
        \draw[ultra thick, black] (3, -8.2) -- (3, -7);
        \draw[ultra thick, black] (3, 8.2) -- (3, 7);
        \draw[ultra thick, green] (3, 1.25) -- (3, 7);
        \draw[ultra thick, green] (3, -1.25) -- (3, -7);
        \draw[ultra thick, blue, fill= white] (3, 0) circle (1.5);
	\draw[ultra thick, black, fill= white] (-1.5,0) ellipse (2.6 and 1) node{$a, b-1$};
	\draw[ultra thick, black, fill= white] (0,-7) ellipse (4 and 1) node{$a, b$};
	\draw[ultra thick, black, fill= white] (0,7) ellipse (4 and 1) node{$a,b$};
\end{tikzpicture}
\end{equation}
becomes
\begin{equation}
\begin{tikzpicture}[scale=.3, anchorbase]
        \draw[ultra thick, green] (3, 1.25) -- (3, 7);
        \draw[ultra thick, green] (3, -1.25) -- (3, -7);
        \draw[ultra thick, blue, fill= white] (3, 0) circle (1.5);
	\draw[ultra thick, black, fill= white] (0,-7) ellipse (4 and 1) ;
	\draw[ultra thick, black, fill= white] (0,7) ellipse (4 and 1) ;
\end{tikzpicture}
\end{equation}

We define $\kappa_{\lambda, \varpi_a}= 1$ for all $\lambda\in X_+$ and $a\in \lbrace \blues, \greent\rbrace$. Whenever $\mu \notin S_{(a, b), \varpi}$, we set $\kappa_{(a,b), \mu}^{-1}$ equal to zero. This results in the following initial conditions for our recursion:
\begin{equation}\label{ic1}
\kappa_{(a, b), (1, -1)}^{-1} = 0 \ \ \ \ \ \text{when} \ \ \ \ \ b= 0,
\end{equation}
\begin{equation}\label{ic2}
\kappa_{(a, b), (-1, 0)}^{-1} = \kappa_{(a,b), (-1, 1)}^{-1}= 0 \ \ \ \ \ \text{when} \ \ \ \ \ a= 0, 
\end{equation}
\begin{equation}\label{ic3}
\kappa_{(a, b), (0,0)}^{-1} = 0 \ \ \ \ \ \text{when} \ \ \ \ \ a= 0,
\end{equation}
\begin{equation}\label{ic4}
\kappa_{(a,b), (-2, 1)}^{-1} = 0 \ \ \ \ \ \text{when} \ \ \ \ \  a= 0 \ \text{or} \ 1,
\end{equation}
and
\begin{equation}\label{ic5}
\kappa_{(a,b), (0, -1)}^{-1} = \kappa_{(a,b), (2, -1)}^{-1} = 0 \ \ \ \ \ \text{when} \ \ \ \ \ b= 0. 
\end{equation}
\end{notation}

\begin{prop}\label{recursionprop}
The $\kappa_{\lambda, \mu}$'s satisfy the following relations.
\begin{equation}
\kappa_{(a, b), (1, 0)} = 1
\end{equation}
\begin{equation}
\kappa_{(a, b), (0, 1)} = 1
\end{equation}
\begin{equation}\label{rr1}
\kappa_{(a,b), (-1, 1)} = -[2] -\dfrac{1}{\kappa_{(a-1, b), (-1,1)}}
\end{equation}
\begin{equation}\label{rr2}
\kappa_{(a, b), (2, -1)} = -\dfrac{[4]}{[2]} - \dfrac{1}{\kappa_{(a, b-1), (2, -1)}}
\end{equation}
\begin{equation}\label{rr3}
\kappa_{(a,b), (0,0)} = \dfrac{[5]}{[2]} - \dfrac{\kappa_{(a-2, b+1), (2, -1)}}{\kappa_{(a-1, b), (-1, 1)}^{-1}} - \dfrac{1}{\kappa_{(a-1, b), (1, -1)}}
\end{equation}
\begin{equation}\label{rr4}
\kappa_{(a,b), (1, -1)} = \dfrac{[5]}{[2]} - \dfrac{\kappa_{(a+2, b-2), (-1, 1)}}{ \kappa_{(a, b-1), (2,-1)}}- \dfrac{1}{[2]^2\kappa_{(a, b-1), (0,0)}}
\end{equation}
\begin{equation}\label{rr5}
\begin{split}
\kappa_{(a,b), (-2,1)} = \dfrac{[5]}{[2]} \kappa_{(a-1, b), (-1, 1)} &- (-[2]- \dfrac{1}{\kappa_{(a-2,b), (-1,1)}})\dfrac{\kappa_{(a-1,b),(-1, 1)}}{\kappa_{(a-1, b),(-1,0)}} \\
&- \dfrac{\kappa_{(a-2, b+1), (0,0)}}{\kappa_{(a-2, b), (-1, 1)}^2 \kappa_{(a-1, b), (-1, 1)}}
\end{split}
\end{equation}
\begin{equation}\label{rr6}
\kappa_{(a,b),(-1,0)} = -\dfrac{[6][2]}{[3]}- \dfrac{1}{\kappa_{(a-1, b), (-1,0)}} - \dfrac{\kappa_{(a-2, b+1), (1, -1)}}{\kappa_{(a-1, b), (-1, 1)}} - \dfrac{\kappa_{(a, b-1), (-1, 1)}}{\kappa_{(a-1, b), (1, -1)}}
\end{equation}
\begin{equation}\label{rr7}
\kappa_{(a, b), (0,-1)} = \dfrac{[6][5]}{[3][2]} - \dfrac{1}{\kappa_{(a,b-1), (0, -1)}} - \dfrac{\kappa_{(a+2, b-2), (-2, 1)}}{ \kappa_{(a, b-1), (2, -1)}} - \dfrac{\kappa_{(a, b-1), (0,0)}}{\kappa_{(a, b-1), (0,0)}} 
- \dfrac{\kappa_{(a-2, b), (2, -1)}}{\kappa_{(a, b-1), (-2, 1)}}
\end{equation}
\end{prop}

\begin{proof} We will use the established properties of clasps to derive the recursion relations. Recall that $\kappa_{\lambda, \mu}$ is the coefficient of the neutral map in $K_{\lambda, \mu}$. 

\begin{equation}
K_{\lambda, \mu}= \quad
\begin{tikzpicture}[scale=.3, anchorbase]
	\draw[ultra thick, black] (-3, 0) -- (-3,7);
	\draw[ultra thick, black] (3, -3.8) -- (3, 3.8);
	\draw[ultra thick, black] (0, -3.8) -- (0, 3.8);
        \draw[ultra thick, black] (-3,0) -- (-3, -7);
	\draw[ultra thick, black] (1.5, -3.8) -- (1.5, -7);
	\draw[ultra thick, black] (1.5, 3.8) -- (1.5, 7);
	\draw[ultra thick, black, fill= white] (0,7) ellipse (4 and 1) node{$\lambda + \mu$};
	\draw[black, ultra thick, fill=white] (-1.5,2.8) -- (4.5,2.8) -- (3.5,4.8) -- (-.5, 4.8) -- cycle;
        \node at (1.5, 3.8) {$L_{\mu}$};
	\draw[ultra thick, black, fill= white] (-1.5,0) ellipse (3 and 1) node{$\lambda$};
	\draw[ultra thick, black, fill= white] (0,-7) ellipse (4 and 1) node{$\lambda + \mu$};
	\draw[black, ultra thick, fill=white] (-1.5,-2.8) -- (4.5,-2.8) -- (3.5,-4.8) -- (-.5, -4.8) -- cycle;
        \node at (1.5, -3.8) {$\mathbb{D}(L_{\mu})$};
\end{tikzpicture}
\end{equation}
\noindent Using the equation 
\begin{equation}
\lambda = (\lambda - \varpi)\ot \varpi - \sum_{\nu\in \wt V(\varpi)- \lbrace \varpi\rbrace} \kappa_{\lambda- \varpi, \nu}^{-1} \mathbb{D} E_{\lambda- \varpi, \nu} \circ E_{\lambda- \varpi, \nu}
\end{equation}
\noindent to rewrite the $\lambda$ clasp, and then using clasp absorption \eqref{claspabs}, we can rewrite $K_{\lambda, \mu}$ as

\begin{equation}\label{LIFclaspreplace}
K_{\lambda, \mu}= \quad
\begin{tikzpicture}[scale=.3, anchorbase]
	\draw[ultra thick, black] (-3, 0) -- (-3,7.7);
	\draw[ultra thick, black] (3, -4.5) -- (3, 4.5);
	\draw[ultra thick, black] (0, -4.5) -- (0, 4.5);
        \draw[ultra thick, black] (-3,0) -- (-3, -7.7);
	\draw[ultra thick, black] (1.5, -4.5) -- (1.5, -7.7);
	\draw[ultra thick, black] (1.5, 4.5) -- (1.5, 7.7);
	\draw[ultra thick, black] (1, 4.5) -- (1, -4.5);
	\draw[ultra thick, black, fill= white] (0,7.7) ellipse (4 and 1) node{$\lambda + \mu$};
	\draw[black, ultra thick, fill=white] (-1.5,3.9) -- (4.5,3.9) -- (3.5,5.9) -- (-.5, 5.9) -- cycle;
        \node at (1.5, 4.9) {$L_{\mu}$};
	\draw[ultra thick, black, fill= white] (-2.5,0) ellipse (3 and .9) node{$\lambda-\varpi$};
	\draw[ultra thick, black, fill= white] (0,-7.7) ellipse (4 and 1) node{$\lambda + \mu$};
	\draw[black, ultra thick, fill=white] (-1.5,-3.9) -- (4.5,-3.9) -- (3.5,-5.9) -- (-.5, -5.9) -- cycle;
        \node at (1.5, -4.9) {$\mathbb{D}(L_{\mu})$};
\end{tikzpicture}
\quad - \sum_{\nu\in \wt V(\varpi)- \lbrace \varpi\rbrace}\kappa_{\lambda-\varpi, \nu}^{-1} \quad
\begin{tikzpicture}[scale=.3, anchorbase]
	\draw[ultra thick, black] (-3, 0) -- (-3,7.7);
	\draw[ultra thick, black] (3, -4.5) -- (3, 4.5);
	\draw[ultra thick, black] (0, -4.5) -- (0, 4.5);
        \draw[ultra thick, black] (-3,0) -- (-3, -7.7);
	\draw[ultra thick, black] (1.5, -4.5) -- (1.5, -7.7);
	\draw[ultra thick, black] (1.5, 4.5) -- (1.5, 7.7);
	\draw[ultra thick, black] (1, 4.5) -- (1, 2.5);
	\draw[ultra thick, black] (1, -2.5) -- (1, -4.5);
	\draw[ultra thick, black, fill= white] (0,7.7) ellipse (4 and 1) node{$\lambda + \mu$};
	\draw[black, ultra thick, fill=white] (-1.5,3.9) -- (4.5,3.9) -- (3.5,5.9) -- (-.5, 5.9) -- cycle;
        \node at (1.5, 4.9) {$L_{\mu}$};
	\draw[ultra thick, black, fill= white] (-2.5,0) ellipse (3 and .9) node{$\lambda-\varpi$};
	\draw[ultra thick, black, fill= white] (0,-7.7) ellipse (4 and 1) node{$\lambda + \mu$};
	\draw[black, ultra thick, fill=white] (-1.5,-3.9) -- (4.5,-3.9) -- (3.5,-5.9) -- (-.5, -5.9) -- cycle;
	\draw[black, ultra thick, fill=white] (-2.7,-3.2) -- (2.7,-3.2) -- (1.7,-1.3) -- (-1.7, -1.3) -- cycle;
	\draw[black, ultra thick, fill=white] (-2.7,3.2) -- (2.7,3.2) -- (1.7,1.3) -- (-1.7, 1.3) -- cycle;
        \node at (1.5, -4.9) {$\mathbb{D}(L_{\mu})$};
        \node at (0, -2.3) {$L_{\nu}$};
        \node at (0, 2.3) {$\mathbb{D}(L_{\nu})$};
\end{tikzpicture}
\end{equation}

Having established the general pattern one follows to derive these recursive formulas, we proceed to apply it for each $\kappa_{(a, b), \mu}$. 

To compute $\kappa_{(a, b), (-1, 1)}$, we resolve the $(a, b)$ clasp in $K_{(a,b), (-1, 1)}$ as in \eqref{LIFclaspreplace}. Since strictly lower diagrams are orthogonal to clasps \eqref{clasporth}, we find
\begin{equation}
\begin{tikzpicture}[scale=.3, anchorbase]
	\draw[ultra thick, blue] (-2, -4) -- (-2,2);
	\draw[ultra thick, green] (2, -4) -- (2,0);
	\draw[ultra thick, blue] (1,1) arc (-180:0:1);
	\draw[ultra thick, blue] (1,2) arc (0:180:1.5);
	\draw[ultra thick, green] (2, -4) -- (2,0);
	\draw[ultra thick, blue] (1, 1) -- (1,2);
	\draw[ultra thick, blue] (3, 1) -- (3,5);
	\draw[ultra thick, black, fill= white] (0,-3) ellipse (4 and 1) node{};
\end{tikzpicture}
\quad = 0
\end{equation}
and
\begin{equation}
\begin{tikzpicture}[scale=.3, anchorbase]
	\draw[ultra thick, green] (-2, -4) -- (-2,2);
	\draw[ultra thick, green] (2, -4) -- (2,0);
	\draw[ultra thick, blue] (1,1) arc (-180:0:1);
	\draw[ultra thick, blue] (1,2) arc (0:90:1.5);
	\draw[ultra thick, green] (-.5,3.5) arc (90:180:1.5);
	\draw[ultra thick, blue] (-.5, 3.5) -- (-.5,5);
	\draw[ultra thick, green] (2, -4) -- (2,0);
	\draw[ultra thick, blue] (1, 1) -- (1,2);
	\draw[ultra thick, blue] (3, 1) -- (3,5);
	\draw[ultra thick, black, fill= white] (0,-3) ellipse (4 and 1) node{};
\end{tikzpicture}
\quad = 0.
\end{equation}
So after expanding the $(a,b)$ clasp in $K_{(a, b), (-1, 1)}$ the only $\nu\in \wt V(\blues)$ which contributes to the sum in \eqref{LIFclaspreplace} is $(-1, 1)$. This means we can rewrite $K_{(a, b), (-1, 1)}$ as follows. 
\begin{equation}
\begin{tikzpicture}[scale=.25, anchorbase]
	\draw[ultra thick, green] (1.5, -2) -- (1.5, -7);
	\draw[ultra thick, green] (1.5, 2) -- (1.5, 7);
	\draw[ultra thick, blue] (1.5,0) ellipse (2 and 2) node{};
	\draw[ultra thick, black, fill= white] (0,7) ellipse (4 and 1) node{};
	\draw[ultra thick, black, fill= white] (-1.5,0) ellipse (3 and 1) node{};
	\draw[ultra thick, black, fill= white] (0,-7) ellipse (4 and 1) node{};
\end{tikzpicture}
\quad= \quad
\begin{tikzpicture}[scale=.25, anchorbase]
	\draw[ultra thick, green] (1.5, -2) -- (1.5, -7);
	\draw[ultra thick, green] (1.5, 2) -- (1.5, 7);
	\draw[ultra thick, blue] (1.5,0) ellipse (2 and 2) node{};
	\draw[ultra thick, black, fill= white] (0,7) ellipse (4 and 1) node{};
	\draw[ultra thick, black, fill= white] (0,-7) ellipse (4 and 1) node{};
\end{tikzpicture}
\quad - \kappa_{(a-1, b), (-1, 1)}^{-1} \quad
\begin{tikzpicture}[scale=.25, anchorbase]
	\draw[ultra thick, green] (1.5, -4) -- (1.5, -7);
	\draw[ultra thick, green] (1.5, 4) -- (1.5, 7);
	\draw[ultra thick, blue] (.5,-3) arc (-180:0:1);
	\draw[ultra thick, blue] (-1.5,-3) arc (180:0:1);
	\draw[ultra thick, blue] (.5,3) arc (180:0:1);
	\draw[ultra thick, blue] (-1.5,3) arc (-180:0:1);
	\draw[ultra thick, blue] (2.5, -3) -- (2.5, 3);
	\draw[ultra thick, green] (-.5, -2) -- (-.5, 2);
	\draw[ultra thick, blue] (-1.5, -3) -- (-1.5, -7);
	\draw[ultra thick, blue] (-1.5, 3) -- (-1.5, 7);
	\draw[ultra thick, black, fill= white] (0,7) ellipse (4 and 1) node{};
	\draw[ultra thick, black, fill= white] (-1.5,0) ellipse (3 and 1) node{};
	\draw[ultra thick, black, fill= white] (0,-7) ellipse (4 and 1) node{};
\end{tikzpicture}
\end{equation}
Using neutral ladder absorption \eqref{neutralabs} and clasp absorption \eqref{claspabs}, we deduce $\kappa_{(a,b), (-1, 1)} = -[2] -\kappa_{(a-1, b), (-1,1)}^{-1}$.

We proceed similarly with the remaining $K_{\lambda, \mu}$. It is useful to note that from the $H\equiv I$ relation \eqref{H=I} we have the following "clasped" relation. 
\begin{equation}\label{h=ionclasp}
\begin{tikzpicture}[scale=.25, anchorbase]
	\draw[ultra thick, blue] (2, -3) -- (2, 3);
	\draw[ultra thick, blue] (-2, -3) -- (-2, 3);
	\draw[ultra thick, green] (-2, 0) -- (2, 0);
	\draw[ultra thick, black, fill= white] (0,-3) ellipse (4 and 1) node{};
\end{tikzpicture}
\quad = \dfrac{1}{[2]}\quad
\begin{tikzpicture}[scale=.25, anchorbase]
	\draw[ultra thick, blue] (2, -3) -- (2, 3);
	\draw[ultra thick, blue] (-2, -3) -- (-2, 3);
	\draw[ultra thick, black, fill= white] (0,-3) ellipse (4 and 1) node{};
\end{tikzpicture}
\end{equation}

Similarly, using the $H\equiv I$ relation \eqref{H=I}, orthogonality to clasps \eqref{clasporth}, and neutral map absorption \eqref{neutralabs} we can also deduce the following. 
\begin{equation}\label{offbrandneutral}
\begin{tikzpicture}[scale=.3, anchorbase]
	\draw[ultra thick, blue] (2.5, -4) -- (2.5, 0);
	\draw[ultra thick, blue] (0, -4) -- (0, 0);
	\draw[ultra thick, blue] (2.5, 0) -- (0, 0);
	\draw[ultra thick, green] (2.5, 0) -- (2.5, 4);
	\draw[ultra thick, green] (-2.5, -4) -- (-2.5, 2);
	\draw[ultra thick, green] (0, 0) -- (0, 2);
	\draw[ultra thick, blue] (-2.5, 2) -- (-2.5, 4);
	\draw[ultra thick, blue] (0, 2) -- (0, 4);
	\draw[ultra thick, blue] (-2.5, 2) -- (0, 2);
	\draw[ultra thick, black, fill= white] (0,-3) ellipse (4 and 1) node{};
\end{tikzpicture}
\quad = \quad
\begin{tikzpicture}[scale=.3, anchorbase]
	\draw[ultra thick, green] (-2.5, -4) -- (-2.5, 0);
	\draw[ultra thick, blue] (-2.5, 0) -- (-2.5, 4);
	\draw[ultra thick, blue] (-2.5, 0) -- (0, 0);
	\draw[ultra thick, blue] (0, -4) -- (0, 0);
	\draw[ultra thick, green] (0, 0) -- (0, 2);
	\draw[ultra thick, blue] (0, 2) -- (0, 4);
	\draw[ultra thick, blue] (0, 2) -- (2.5, 2);
	\draw[ultra thick, blue] (2.5, -4) -- (2.5, 2);
	\draw[ultra thick, green] (2.5, 2) -- (2.5, 4);
	\draw[ultra thick, black, fill= white] (0,-3) ellipse (4 and 1) node{};
\end{tikzpicture}
\quad = \quad 
\begin{tikzpicture}[scale=.3, anchorbase]
	\draw[ultra thick, green] (-2.5, -4) -- (-2.5, 4);
	\draw[ultra thick, blue] (0, -4) -- (0, 4);
	\draw[ultra thick, blue] (2.5, -4) -- (2.5, 4);
	\draw[ultra thick, black, fill= white] (0,-3) ellipse (4 and 1) node{};
\end{tikzpicture}
\end{equation}

For $\kappa_{(a,b), (2, -1)}$, we begin by observing that by clasp orthogonality \eqref{clasporth}
\begin{equation}
\begin{tikzpicture}[scale=.3, anchorbase]
	\draw[ultra thick, green] (-3, -3.5) -- (-3,2);
	\draw[ultra thick, green] (0,2) arc (0:180:1.5);
	\draw[ultra thick, blue] (0, -3.5) -- (0,2);
	\draw[ultra thick, blue] (0, 2) -- (3,2);
	\draw[ultra thick, blue] (3, -3.5) -- (3,2);
	\draw[ultra thick, green] (3, 2) -- (3,5);
	\draw[ultra thick, black, fill= white] (0,-3) ellipse (4 and 1) node{};
\end{tikzpicture}
\quad = 0, \quad
\end{equation}
and by \eqref{h=ionclasp} and clasp orthogonality \eqref{clasporth},
\begin{equation}
\begin{tikzpicture}[scale=.3, anchorbase]
	\draw[ultra thick, blue] (-3, -3.5) -- (-3,2);
	\draw[ultra thick, green] (0,2) arc (0:90:1.5);
	\draw[ultra thick, blue] (-1.5,3.5) arc (90:180:1.5);
	\draw[ultra thick, blue] (-1.5, 3.5) -- (-1.5,5);
	\draw[ultra thick, blue] (0, -4) -- (0,2);
	\draw[ultra thick, blue] (0, 2) -- (3,2);
	\draw[ultra thick, blue] (3, -3.5) -- (3,2);
	\draw[ultra thick, green] (3, 2) -- (3,5);
	\draw[ultra thick, black, fill= white] (0,-3) ellipse (4 and 1) node{};
\end{tikzpicture}
\quad = 0.
\end{equation}
Using these observations and \eqref{offbrandneutral}, $K_{(a, b), (2, -1)}$ can be resolved as follows. 
\begin{equation}
\begin{tikzpicture}[scale=.25, anchorbase]
	\draw[ultra thick, blue] (0, -3) -- (0, -7);
	\draw[ultra thick, blue] (0, 3) -- (0, 7);
	\draw[ultra thick, blue] (3, -3) -- (3, -7);
	\draw[ultra thick, blue] (3, 3) -- (3, 7);
	\draw[ultra thick, blue] (0, -3) -- (3, -3);
	\draw[ultra thick, blue] (0, 3) -- (3, 3);
	\draw[ultra thick, green] (0, -3) -- (0, 3);
	\draw[ultra thick, green] (3, -3) -- (3, 3);
	\draw[ultra thick, black, fill= white] (0,7) ellipse (4 and 1) node{};
	\draw[ultra thick, black, fill= white] (-1.5,0) ellipse (3 and 1) node{};
	\draw[ultra thick, black, fill= white] (0,-7) ellipse (4 and 1) node{};
\end{tikzpicture}
\quad= \quad
\begin{tikzpicture}[scale=.25, anchorbase]
	\draw[ultra thick, blue] (0, -3) -- (0, -7);
	\draw[ultra thick, blue] (0, 3) -- (0, 7);
	\draw[ultra thick, blue] (3, -3) -- (3, -7);
	\draw[ultra thick, blue] (3, 3) -- (3, 7);
	\draw[ultra thick, blue] (0, -3) -- (3, -3);
	\draw[ultra thick, blue] (0, 3) -- (3, 3);
	\draw[ultra thick, green] (0, -3) -- (0, 3);
	\draw[ultra thick, green] (3, -3) -- (3, 3);
	\draw[ultra thick, black, fill= white] (0,7) ellipse (4 and 1) node{};
	\draw[ultra thick, black, fill= white] (0,-7) ellipse (4 and 1) node{};
\end{tikzpicture}
\quad - \kappa_{(a, b-1), (2, -1)}^{-1} \quad
\begin{tikzpicture}[scale=.25, anchorbase]
	\draw[ultra thick, blue] (0, -7) -- (0, 7);
	\draw[ultra thick, blue] (3, -7) -- (3, 7);
	\draw[ultra thick, green] (-3, -7) -- (-3, 7);
	\draw[ultra thick, black, fill= white] (0,7) ellipse (4 and 1) node{};
	\draw[ultra thick, black, fill= white] (0,-7) ellipse (4 and 1) node{};
\end{tikzpicture}
\end{equation}
Then we apply the $H \equiv I$ relation \eqref{H=I} to find
\begin{equation}
\begin{tikzpicture}[scale=.25, anchorbase]
	\draw[ultra thick, blue] (0, -3) -- (0, -7);
	\draw[ultra thick, blue] (0, 3) -- (0, 7);
	\draw[ultra thick, blue] (3, -3) -- (3, -7);
	\draw[ultra thick, blue] (3, 3) -- (3, 7);
	\draw[ultra thick, blue] (0, -3) -- (3, -3);
	\draw[ultra thick, blue] (0, 3) -- (3, 3);
	\draw[ultra thick, green] (0, -3) -- (0, 3);
	\draw[ultra thick, green] (3, -3) -- (3, 3);
	\draw[ultra thick, black, fill= white] (0,7) ellipse (4 and 1) node{};
	\draw[ultra thick, black, fill= white] (0,-7) ellipse (4 and 1) node{};
\end{tikzpicture}
\quad = \quad
\begin{tikzpicture}[scale=.25, anchorbase]
	\draw[ultra thick, green] (0, 0) -- (1.5, 0);
	\draw[ultra thick, blue] (0, -7) -- (0, 7);
	\draw[ultra thick, blue] (3, -3) -- (3, -7);
	\draw[ultra thick, blue] (3, 3) -- (3, 7);
	\draw[ultra thick, green] (3, -3) -- (3, 3);
	\draw[ultra thick, blue] (3,-3) .. controls (1,0) .. (3,3);
	\draw[ultra thick, black, fill= white] (0,7) ellipse (4 and 1) node{};
	\draw[ultra thick, black, fill= white] (0,-7) ellipse (4 and 1) node{};
\end{tikzpicture}
\quad - \dfrac{1}{[2]}\quad
\begin{tikzpicture}[scale=.25, anchorbase]
	\draw[ultra thick, blue] (0, -7) -- (0, 7);
	\draw[ultra thick, blue] (3, -3) -- (3, -7);
	\draw[ultra thick, blue] (3, 3) -- (3, 7);
	\draw[ultra thick, green] (3, -3) -- (3, 3);
	\draw[ultra thick, blue] (3,-3) .. controls (1,0) .. (3,3);
	\draw[ultra thick, black, fill= white] (0,7) ellipse (4 and 1) node{};
	\draw[ultra thick, black, fill= white] (0,-7) ellipse (4 and 1) node{};
\end{tikzpicture}
\end{equation}
Here the vanishing of the third term is due to clasp orthogonality \eqref{clasporth}. If we apply the $H \equiv I$ relation \eqref{H=I} again, then by the monogon relation \eqref{monogon} and clasp orthogonality \eqref{clasporth} we can rewrite the right hand side as follows. 
\begin{equation}\label{ladderrhs}
\begin{tikzpicture}[scale=.25, anchorbase]
	\draw[ultra thick, green] (-1, 0) -- (0, 0);
	\draw[ultra thick, green] (3, 0) -- (2.25, 0);
	\draw[ultra thick, blue] (-1, -7) -- (-1, 7);
	\draw[ultra thick, blue] (3, 7) -- (3, -7);
	\draw[ultra thick, blue] (1,0) ellipse (1 and 1) node{};
	\draw[ultra thick, black, fill= white] (0,7) ellipse (4 and 1) node{};
	\draw[ultra thick, black, fill= white] (0,-7) ellipse (4 and 1) node{};
\end{tikzpicture}
\quad +\dfrac{1}{[2]} \quad
\begin{tikzpicture}[scale=.25, anchorbase]
	\draw[ultra thick, blue] (-1, -7) -- (-1, 7);
	\draw[ultra thick, blue] (3, 7) -- (3, -7);
	\draw[ultra thick, green] (-1, 0) -- (3, 0);
	\draw[ultra thick, black, fill= white] (0,7) ellipse (4 and 1) node{};
	\draw[ultra thick, black, fill= white] (0,-7) ellipse (4 and 1) node{};
\end{tikzpicture}
\quad - \dfrac{[5]}{[2]^2}\quad
\begin{tikzpicture}[scale=.25, anchorbase]
	\draw[ultra thick, blue] (-1, -7) -- (-1, 7);
	\draw[ultra thick, blue] (3, 7) -- (3, -7);
	\draw[ultra thick, black, fill= white] (0,7) ellipse (4 and 1) node{};
	\draw[ultra thick, black, fill= white] (0,-7) ellipse (4 and 1) node{};
\end{tikzpicture}
\end{equation}
Finally, using the bigon relation and the clasped $H \equiv I$ relation \eqref{H=I} we can rewrite \eqref{ladderrhs} as 
\begin{equation}
-[2] \dfrac{1}{[2]} + \dfrac{1}{[2]^2}-\dfrac{[5]}{[2]^2}= -\dfrac{[4]}{[2]}
\end{equation}
times the clasp, and then conclude that
\begin{equation}
\kappa_{(a, b), (2, -1)} = -\dfrac{[4]}{[2]} - \kappa_{(a, b-1), (2, -1)}^{-1}.
\end{equation}

To compute $\kappa_{(a,b), (0,0)}$ we will expand the middle clasp in $K_{(a, b), (0,0)}$.
\begin{equation}
K_{(a, b), (0, 0)}= \quad
\begin{tikzpicture}[scale=.25, anchorbase]
	\draw[ultra thick, blue] (1.5, -2) -- (1.5, -7);
	\draw[ultra thick, blue] (1.5, 2) -- (1.5, 7);
	\draw[ultra thick, blue] (1.5,-2) arc (270: 90: 2);
	\draw[ultra thick, green] (1.5,-2) arc (-90: 90: 2);
	\draw[ultra thick, black, fill= white] (0,7) ellipse (4 and 1) node{};
	\draw[ultra thick, black, fill= white] (-1.5,0) ellipse (3 and 1) node{};
	\draw[ultra thick, black, fill= white] (0,-7) ellipse (4 and 1) node{};
\end{tikzpicture}
\end{equation}
Since
\begin{equation}
\begin{tikzpicture}[scale=.25, anchorbase]
	\draw[ultra thick, blue] (-2, -4) -- (-2,2);
	\draw[ultra thick, blue] (2, -4) -- (2,0);
	\draw[ultra thick, blue] (1,1) arc (-180:-90:1);
	\draw[ultra thick, green] (2, 0) arc (-90: 0:1);
	\draw[ultra thick, blue] (1,2) arc (0:180:1.5);
	\draw[ultra thick, blue] (1, 1) -- (1,2);
	\draw[ultra thick, green] (3, 1) -- (3,5);
	\draw[ultra thick, black, fill= white] (0,-3) ellipse (4 and 1) node{};
\end{tikzpicture}
\quad = 0,
\end{equation}
we can rewrite $K_{(a, b), (0,0)}$ as follows. 
\begin{equation}\label{kab00expanded}
\begin{tikzpicture}[scale=.25, anchorbase]
	\draw[ultra thick, blue] (1.5, -2) -- (1.5, -7);
	\draw[ultra thick, blue] (1.5, 2) -- (1.5, 7);
	\draw[ultra thick, blue] (1.5,-2) arc (270: 90: 2);
	\draw[ultra thick, green] (1.5,-2) arc (-90: 90: 2);
	\draw[ultra thick, black, fill= white] (0,7) ellipse (4 and 1) node{};
	\draw[ultra thick, black, fill= white] (0,-7) ellipse (4 and 1) node{};
\end{tikzpicture}
\quad - \kappa_{(a-1, b), (-1, 1)}^{-1} \quad
\begin{tikzpicture}[scale=.25, anchorbase]
	\draw[ultra thick, green] (1.5,-4) arc (-90:0:1);
	\draw[ultra thick, blue] (0.5,-3) arc (-180:-90:1);
	\draw[ultra thick, blue] (-1.5,-3) arc (180:0:1);
	\draw[ultra thick, blue] (.5,3) arc (180:90:1);
	\draw[ultra thick, green] (2.5,3) arc (0:90:1);
	\draw[ultra thick, blue] (-1.5,3) arc (-180:0:1);
	\draw[ultra thick, blue] (1.5, -4) -- (1.5, -7);
	\draw[ultra thick, blue] (1.5, 4) -- (1.5, 7);
	\draw[ultra thick, green] (2.5, -3) -- (2.5, 3);
	\draw[ultra thick, green] (-.5, -2) -- (-.5, 2);
	\draw[ultra thick, blue] (-1.5, -3) -- (-1.5, -7);
	\draw[ultra thick, blue] (-1.5, 3) -- (-1.5, 7);
	\draw[ultra thick, black, fill= white] (0,7) ellipse (4 and 1) node{};
	\draw[ultra thick, black, fill= white] (-1.5,0) ellipse (3 and 1) node{};
	\draw[ultra thick, black, fill= white] (0,-7) ellipse (4 and 1) node{};
\end{tikzpicture}
\quad - \kappa_{(a-1, b), (1, -1)}^{-1} \quad
\begin{tikzpicture}[scale=.25, anchorbase]
	\draw[ultra thick, green] (1.5,-4) arc (-90:0:1);
	\draw[ultra thick, blue] (0.5,-3) arc (-180:-90:1);
	\draw[ultra thick, green] (-1.5,-3) arc (180:90:1);
	\draw[ultra thick, blue] (.5,-3) arc (0:90:1);
	\draw[ultra thick, blue] (.5,3) arc (0:-90:1);
	\draw[ultra thick, blue] (.5,3) arc (180:90:1);
	\draw[ultra thick, green] (2.5,3) arc (0:90:1);
	\draw[ultra thick, green] (-1.5,3) arc (-180:-90:1);
	\draw[ultra thick, blue] (1.5, -4) -- (1.5, -7);
	\draw[ultra thick, blue] (1.5, 4) -- (1.5, 7);
	\draw[ultra thick, green] (2.5, -3) -- (2.5, 3);
	\draw[ultra thick, blue] (-.5, -2) -- (-.5, 2);
	\draw[ultra thick, green] (-1.5, -3) -- (-1.5, -7);
	\draw[ultra thick, green] (-1.5, 3) -- (-1.5, 7);
	\draw[ultra thick, black, fill= white] (0,7) ellipse (4 and 1) node{};
	\draw[ultra thick, black, fill= white] (-1.5,0) ellipse (3 and 1) node{};
	\draw[ultra thick, black, fill= white] (0,-7) ellipse (4 and 1) node{};
\end{tikzpicture}
\end{equation}
Observing that the second term in \eqref{kab00expanded} is $K_{(a-2, b+1), (2, -1)}$, and then using neutral map absorption \eqref{neutralabs} and clasp absorption \eqref{claspabs} for the third term in \eqref{kab00expanded}, we deduce that
\begin{equation}
\kappa_{(a,b), (0,0)} = \dfrac{[5]}{[2]} - \kappa_{(a-1, b), (-1, 1)}^{-1} \kappa_{(a-2, b+1), (2, -1)} - \kappa_{(a-1, b), (1, -1)}^{-1}. 
\end{equation}

To compute $\kappa_{(a,b),(1, -1)}$ we expand the middle clasp in $K_{(a, b), (1, -1)}$.
\begin{equation}
\begin{tikzpicture}[scale=.25, anchorbase]
	\draw[ultra thick, blue] (1.5, -2) -- (1.5, -7);
	\draw[ultra thick, blue] (1.5, 2) -- (1.5, 7);
	\draw[ultra thick, green] (1.5,-2) arc (270: 90: 2);
	\draw[ultra thick, blue] (1.5,-2) arc (-90: 90: 2);
	\draw[ultra thick, black, fill= white] (0,7) ellipse (4 and 1) node{};
	\draw[ultra thick, black, fill= white] (-1.5,0) ellipse (3 and 1) node{};
	\draw[ultra thick, black, fill= white] (0,-7) ellipse (4 and 1) node{};
\end{tikzpicture}
\end{equation}
We begin by using the $H \equiv I$ relation \eqref{H=I} and clasp orthogonality \eqref{clasporth}, followed by neutral absorption \eqref{neutralabs}, to calculate the following. 
\begin{equation}
\begin{tikzpicture}[scale=.25, anchorbase]
	\draw[ultra thick, blue] (1.5,-4) arc (-90:0:1);
	\draw[ultra thick, green] (0.5,-3) arc (-180:-90:1);
	\draw[ultra thick, blue] (1.5, -4) -- (1.5, -7);
	\draw[ultra thick, blue] (2.5, -3) -- (2.5, 3);
	\draw[ultra thick, green] (-1.5, -3) -- (-1.5, -7);
	\draw[ultra thick, blue] (-1.5, -3) -- (.5, -3);
	\draw[ultra thick, blue] (-1.5, .5) -- (-1.5, -3);
	\draw[ultra thick, blue] (.5, .5) -- (.5, -3);
	\draw[ultra thick, black, fill= white] (-1.5,0) ellipse (3 and 1) node{};
	\draw[ultra thick, black, fill= white] (0,-7) ellipse (4 and 1) node{};
\end{tikzpicture}
\quad = \quad
\begin{tikzpicture}[scale=.25, anchorbase]
	\draw[ultra thick, green] (1.5, -2) -- (1.5, -4);
	\draw[ultra thick, green] (-1.5, -4) -- (-1.5, -7);
	\draw[ultra thick, blue] (1.5,-2) arc (-90:0:1);
	\draw[ultra thick, blue] (0.5,-1) arc (-180:-90:1);
	\draw[ultra thick, blue] (1.5, -7) -- (1.5, -4);
	\draw[ultra thick, blue] (2.5, -1) -- (2.5, 3);
	\draw[ultra thick, blue] (-1.5, .5) -- (-1.5, -4);
	\draw[ultra thick, blue] (-1.5, -4) -- (1.5, -4);
	\draw[ultra thick, black, fill= white] (-1.5,0) ellipse (3 and 1) node{};
	\draw[ultra thick, black, fill= white] (0,-7) ellipse (4 and 1) node{};
\end{tikzpicture}
\quad = \quad
\begin{tikzpicture}[scale=.25, anchorbase]
	\draw[ultra thick, green] (1.5, -2) -- (1.5, -7);
	\draw[ultra thick, blue] (-1.5, .5) -- (-1.5, -7);
	\draw[ultra thick, blue] (1.5,-2) arc (-90:0:1);
	\draw[ultra thick, blue] (0.5,-1) arc (-180:-90:1);
	\draw[ultra thick, blue] (2.5, -1) -- (2.5, 3);
	\draw[ultra thick, black, fill= white] (-1.5,0) ellipse (3 and 1) node{};
	\draw[ultra thick, black, fill= white] (0,-7) ellipse (4 and 1) node{};
\end{tikzpicture}
\end{equation}

Clasp orthogonality \eqref{clasporth} implies
\begin{equation}
\begin{tikzpicture}[scale=.25, anchorbase]
	\draw[ultra thick, green] (-2, -4) -- (-2,2);
	\draw[ultra thick, blue] (2, -4) -- (2,0);
	\draw[ultra thick, green] (1,1) arc (-180:-90:1);
	\draw[ultra thick, blue] (2, 0) arc (-90: 0:1);
	\draw[ultra thick, green] (1,2) arc (0:180:1.5);
	\draw[ultra thick, green] (1, 1) -- (1,2);
	\draw[ultra thick, blue] (3, 1) -- (3,5);
	\draw[ultra thick, black, fill= white] (0,-3) ellipse (4 and 1) node{};
\end{tikzpicture}
\quad = 0, \quad
\end{equation}
and the $H\equiv I$ relation followed by clasp orthogonality \eqref{clasporth} implies that
\begin{equation}
\begin{tikzpicture}[scale=.3, anchorbase]
	\draw[ultra thick, blue] (-3, -3.5) -- (-3,2);
	\draw[ultra thick, blue] (0,2) arc (0:90:1.5);
	\draw[ultra thick, blue] (-1.5,3.5) arc (90:180:1.5);
	\draw[ultra thick, green] (-1.5, 3.5) -- (-1.5,5);
	\draw[ultra thick, blue] (0, -4) -- (0,2);
	\draw[ultra thick, green] (0, 2) -- (3,2);
	\draw[ultra thick, blue] (3, -3.5) -- (3,2);
	\draw[ultra thick, blue] (3, 2) -- (3,5);
	\draw[ultra thick, black, fill= white] (0,-3) ellipse (4 and 1) node{};
\end{tikzpicture}
\quad = 0.
\end{equation}
Therefore, we can rewrite $K_{(a, b), (1, -1)}$ as follows. 
\begin{equation}\label{kab1-1expanded}
\begin{tikzpicture}[scale=.25, anchorbase]
	\draw[ultra thick, blue] (1.5, -2) -- (1.5, -7);
	\draw[ultra thick, blue] (1.5, 2) -- (1.5, 7);
	\draw[ultra thick, green] (1.5,-2) arc (270: 90: 2);
	\draw[ultra thick, blue] (1.5,-2) arc (-90: 90: 2);
	\draw[ultra thick, black, fill= white] (0,7) ellipse (4 and 1) node{};
	\draw[ultra thick, black, fill= white] (0,-7) ellipse (4 and 1) node{};
\end{tikzpicture}
\quad - \kappa_{(a, b-1), (2, -1)}^{-1} \quad
\begin{tikzpicture}[scale=.25, anchorbase]
	\draw[ultra thick, green] (1.5, -2) -- (1.5, -7);
	\draw[ultra thick, green] (1.5, 2) -- (1.5, 7);
	\draw[ultra thick, blue] (1.5,-2) arc (270: 90: 2);
	\draw[ultra thick, blue] (1.5,-2) arc (-90: 90: 2);
	\draw[ultra thick, black, fill= white] (0,7) ellipse (4 and 1) node{};
	\draw[ultra thick, black, fill= white] (-1.5,0) ellipse (3 and 1) node{};
	\draw[ultra thick, black, fill= white] (0,-7) ellipse (4 and 1) node{};
\end{tikzpicture}
\quad - \kappa_{(a, b-1), (0, 0)}^{-1} \dfrac{1}{[2]^2} \quad
\begin{tikzpicture}[scale=.25, anchorbase]
	\draw[ultra thick, blue] (2.5, -7) -- (2.5, 7);
	\draw[ultra thick, blue] (.5, -7) -- (.5, 7);
	\draw[ultra thick, black, fill= white] (0,7) ellipse (4 and 1) node{};
	\draw[ultra thick, black, fill= white] (0,-7) ellipse (4 and 1) node{};
\end{tikzpicture}
\end{equation}
Identifying the second term in \eqref{kab1-1expanded} as $K_{(a+2, b-2), (-1, 1)}$, we deduce that
\begin{equation}
\kappa_{(a,b), (1, -1)} = \dfrac{[5]}{[2]} - \kappa_{(a, b-1), (2,-1)}^{-1} \kappa_{(a+2, b-2), (-1, 1)} - \dfrac{1}{[2]^2}\kappa_{(a, b-1), (0,0)}^{-1}. 
\end{equation}

\begin{rmk}
Note that at this point we could start solving these recursive relations, as the local intersection forms for the weights $(-1, 1)$ and $(2, -1)$ are linked only to themselves in their recursion relation. While the local intersection forms for the weights $(0,0)$ and $(1, -1)$ have recursions which link them to themselves, each other, and the local intersection forms for the weights $(-1, 1)$ and $(2, -1)$. 
\end{rmk}

Continuing with our derivation of recursive relations for local intersection forms, we expand the middle clasp in \begin{equation}
K_{(a, b), (-1, 0)}= \quad
\begin{tikzpicture}[scale=.25, anchorbase]
	\draw[ultra thick, blue] (1.5,-2) arc (270: 90: 2);
	\draw[ultra thick, blue] (1.5,-2) arc (-90: 90: 2);
	\draw[ultra thick, black, fill= white] (0,7) ellipse (4 and 1) node{};
	\draw[ultra thick, black, fill= white] (-1.5,0) ellipse (3 and 1) node{};
	\draw[ultra thick, black, fill= white] (0,-7) ellipse (4 and 1) node{};
\end{tikzpicture}
\end{equation}
and apply clasp absorption \eqref{claspabs} to deduce
\begin{equation}\label{kab-10expanded}
\begin{split}
\quad K_{(a, b), (-1, 0)}&= \quad
\begin{tikzpicture}[scale=.45, anchorbase]
	\draw[ultra thick, blue] (3,0) ellipse (1 and 1) node{};
	\draw[ultra thick, black, fill= white] (-.5,0) ellipse (2 and .7) node{};
\end{tikzpicture}
\quad - \kappa_{(a-1, b), (-1, 1)}^{-1}\quad
\begin{tikzpicture}[scale=.45, anchorbase]
	\draw[ultra thick, blue] (1,3.5) -- (1, 2.4);
	\draw[ultra thick, blue] (1,-3.5) -- (1, -2.4);
	\draw[ultra thick, blue] (2,2.4) arc (180:0:.5);
	\draw[ultra thick, blue] (2,2.4) arc (0:-90:.5);
	\draw[ultra thick, blue] (1,2.4) arc (-180:-90:.5);
	\draw[ultra thick, blue] (2,-2.4) arc (0:90:.5);
	\draw[ultra thick, blue] (2,-2.4) arc (-180:0:.5);
	\draw[ultra thick, blue] (1,-2.4) arc (180:90:.5);
	\draw[ultra thick, green] (1.5, -1.9) -- (1.5, 1.9);
	\draw[ultra thick, blue] (3,-2.4) -- (3, 2.4);
	\draw[ultra thick, black, fill= white] (-.5,-3.2) ellipse (2.2 and .7) node{};
	\draw[ultra thick, black, fill= white] (-.5,3.2) ellipse (2.2 and .7) node{};
	\draw[ultra thick, black, fill= white] (0,0) ellipse (2.2 and .7) node{};
\end{tikzpicture}\\
\quad  &- \kappa_{(a-1, b), (1, -1)}^{-1} \quad \begin{tikzpicture}[scale=.45, anchorbase]
	\draw[ultra thick, green] (1,3.5) -- (1, 2.4);
	\draw[ultra thick, green] (1,-3.5) -- (1, -2.4);
	\draw[ultra thick, blue] (2,2.4) arc (180:0:.5);
	\draw[ultra thick, blue] (2, 2.4) arc (0:-90:.5);
	\draw[ultra thick, green] (1,2.4) arc (-180:-90:.5);
	\draw[ultra thick, blue] (2,-2.4) arc (0:90:.5);
	\draw[ultra thick, green] (1,-2.4) arc (180:90:.5);
	\draw[ultra thick, blue] (2,-2.4) arc (-180:0:.5);
	\draw[ultra thick, blue] (1.5, -1.9) -- (1.5, 1.9);
	\draw[ultra thick, blue] (3,-2.4) -- (3, 2.4);
	\draw[ultra thick, black, fill= white] (-.5,-3.2) ellipse (2.2 and .7) node{};
	\draw[ultra thick, black, fill= white] (-.5,3.2) ellipse (2.2 and .7) node{};
	\draw[ultra thick, black, fill= white] (0,0) ellipse (2.2 and .7) node{};
\end{tikzpicture}
\quad - \kappa_{(a-1, b), (-1, 0)}^{-1}\quad \begin{tikzpicture}[scale=.45, anchorbase]
	\draw[ultra thick, blue] (1,3.5) -- (1, 2.4);
	\draw[ultra thick, blue] (1,-3.5) -- (1, -2.4);
	\draw[ultra thick, blue] (2,-2.4) arc (-180:0:.5);
	\draw[ultra thick, blue] (2,2.4) arc (180:0:.5);
	\draw[ultra thick, blue] (2, 2.4) arc (0:-90:.5);
	\draw[ultra thick, blue] (1,2.4) arc (-180:-90:.5);
	\draw[ultra thick, blue] (2,-2.4) arc (0:90:.5);
	\draw[ultra thick, blue] (1,-2.4) arc (180:90:.5);
	\draw[ultra thick, blue] (3,-2.4) -- (3, 2.4);
	\draw[ultra thick, black, fill= white] (-.5,-3.2) ellipse (2.2 and .7) node{};
	\draw[ultra thick, black, fill= white] (-.5,3.2) ellipse (2.2 and .7) node{};
	\draw[ultra thick, black, fill= white] (0,0) ellipse (2.2 and .7) node{};
\end{tikzpicture}
\end{split}
\end{equation}
By identifying the second and third terms on the right hand side of \eqref{kab-10expanded} as $K_{(a-2, b+1), (1, -1)}$ and $K_{(a, b-1), (-1, 1)}$ respectively, we find 
\begin{equation}
\kappa_{(a,b),(-1,0)} = -\dfrac{[6][2]}{[3]} - \dfrac{\kappa_{(a-2, b+1), (1, -1)}}{\kappa_{(a-1, b), (-1, 1)}} - \dfrac{\kappa_{(a, b-1), (-1, 1)}}{\kappa_{(a-1, b), (1, -1)}} - \dfrac{1}{\kappa_{(a-1, b), (-1,0)}}
\end{equation}

\begin{rmk}
Again, we could stop here and solve the recursive relations since the local intersection form for the weight $(-1, 0)$ involves the weights $(-1, 0)$ and the weights we have computed recursions for previously.
\end{rmk}

The antidominant weight in $V(\varpi_2)$ is $(0, -1)$. A calculation similar to the derivation of the recursion for $\kappa_{(a, b), (-1, 0)}$ results in the following
\begin{equation}
\kappa_{(a, b), (0,-1)} = \dfrac{[6][5]}{[3][2]} - - \dfrac{\kappa_{(a-2, b), (2, -1)}}{\kappa_{(a, b-1), (-2, 1)}} - \dfrac{\kappa_{(a, b-1), (0,0)}}{\kappa_{(a, b-1), (0,0)}} 
- \dfrac{\kappa_{(a+2, b-2), (-2, 1)}}{ \kappa_{(a, b-1), (2, -1)}}- \dfrac{1}{\kappa_{(a,b-1), (0, -1)}}. 
\end{equation}

The last local intersection form to resolve is $K_{(a, b), (-2, 1)}$. Recall that 
\begin{equation}
\ak \ot \eval: \DDk\longrightarrow \Fund(U_q^{\ak}(\mathfrak{sp}_4))
\end{equation}
is an equivalence. Also, we know that if a Weyl module $\ak \ot V^{\ak}(\lambda)$ is simple, then we can compute the dimension of homomorphism spaces involving that Weyl module in characteristic zero. Thus, from the calculations in \cite{bodish2020web} we see that
\begin{equation}
\dim \Hom_{\ak \ot \DD}((a-2, b+1), (a, b-1)\ot \greent) = \dim \Hom_{\mathfrak{sp}_4(\mathbb{C})}(V(a-2, b+1), V(a, b-1)\ot V(\varpi_{\greent}))= 0, 
\end{equation}
and it follows that
\begin{equation}\label{claspedtrigon}
\begin{tikzpicture}[scale=.25, anchorbase]
	\draw[ultra thick, blue] (1.5,-4) arc (-90:0:1);
	\draw[ultra thick, blue] (0.5,-3) arc (-180:-90:1);
	\draw[ultra thick, green] (1.5, -4) -- (1.5, -7);
	\draw[ultra thick, blue] (2.5, -3) -- (2.5, 0);
	\draw[ultra thick, blue] (.5, 0) -- (2.5, 0);
	\draw[ultra thick, green] (.5, -3) -- (.5, 0);
	\draw[ultra thick, blue] (.5, 3) -- (.5, 0);
	\draw[ultra thick, green] (2.5, 3) -- (2.5, 0);
	\draw[ultra thick, black, fill= white] (-1.5,-2.5) ellipse (2.5 and .7) node{};
	\draw[ultra thick, black, fill= white] (-1.5,2.5) ellipse (2.5 and .7) node{};
	\draw[ultra thick, black, fill= white] (0,-7) ellipse (4 and 1) node{};
\end{tikzpicture}
\quad = 0.
\end{equation}

When we expand the $(a,b)$ clasp in $K_{(a, b), (-2, 1)}$, one of the terms has \eqref{claspedtrigon} as a sub-diagram and therefore is zero, so we get the following three terms.
\begin{equation}\label{K-21}
\begin{tikzpicture}[scale=.25, anchorbase]
	\draw[ultra thick, green] (1.5, -4) -- (1.5, -7);
	\draw[ultra thick, blue] (1.5, -2) -- (1.5, -4);
	\draw[ultra thick, green] (1.5, 2) -- (1.5, 7);
	\draw[ultra thick, blue] (1.5, 2) -- (1.5, 4);
	\draw[ultra thick, blue] (1.5, -4) ..controls (-3,0) .. (1.5, 4);
	\draw[ultra thick, blue] (1.5,-2) arc (270: 90: 2);
	\draw[ultra thick, green] (1.5,-2) arc (-90: 90: 2);
	\draw[ultra thick, black, fill= white] (0,7) ellipse (4 and 1) node{};
	\draw[ultra thick, black, fill= white] (-3.2,0) ellipse (2.2 and .7) node{};
	\draw[ultra thick, black, fill= white] (0,-7) ellipse (4 and 1) node{};
\end{tikzpicture}
\quad - \kappa_{(a-1, b), (-1, 0)}^{-1} \quad
\begin{tikzpicture}[scale=.25, anchorbase]
	\draw[ultra thick, green] (1.5, -4.5) -- (1.5, -7);
	\draw[ultra thick, green] (1.5, 4.5) -- (1.5, 7);
	\draw[ultra thick, green] (1.5, 1.5) -- (1.5, -1.5);
	\draw[ultra thick, blue] (1.5,-4.5) arc (270: 90: 1.5);
	\draw[ultra thick, blue] (1.5,-4.5) arc (-90: 90: 1.5);
	\draw[ultra thick, blue] (1.5,1.5) arc (270: 90: 1.5);
	\draw[ultra thick, blue] (1.5,1.5) arc (-90: 90: 1.5);
	\draw[ultra thick, black, fill= white] (0,7) ellipse (4 and 1) node{};
	\draw[ultra thick, black, fill= white] (-1,-3) ellipse (2.2 and .7) node{};
	\draw[ultra thick, black, fill= white] (-1,3) ellipse (2.2 and .7) node{};
	\draw[ultra thick, black, fill= white] (0,-7) ellipse (4 and 1) node{};
\end{tikzpicture}
\quad - \kappa_{(a-1, b), (-1, 1)}^{-1}
\begin{tikzpicture}[scale=.25, anchorbase]
	\draw[ultra thick, green] (1.5, -4.5) -- (1.5, -7);
	\draw[ultra thick, green] (1.5, 4.5) -- (1.5, 7);
	\draw[ultra thick, green] (.5, -1.9) ..controls (-.5, 0) .. (.5, 1.9);
	\draw[ultra thick, green] (2.5, -1.9) ..controls (3.5, 0) .. (2.5, 1.9);
	\draw[ultra thick, blue] (1.5,-4.5) arc (270: 90: 1.5);
	\draw[ultra thick, blue] (1.5,-4.5) arc (-90: 90: 1.5);
	\draw[ultra thick, blue] (1.5,1.5) arc (270: 90: 1.5);
	\draw[ultra thick, blue] (1.5,1.5) arc (-90: 90: 1.5);
	\draw[ultra thick, black, fill= white] (0,7) ellipse (4 and 1) node{};
	\draw[ultra thick, black, fill= white] (-1,-3) ellipse (2.2 and .7) node{};
	\draw[ultra thick, black, fill= white] (-1,3) ellipse (2.2 and .7) node{};
	\draw[ultra thick, black, fill= white] (-2,0) ellipse (2.2 and .7) node{};
	\draw[ultra thick, black, fill= white] (0,-7) ellipse (4 and 1) node{};
\end{tikzpicture}
\end{equation}
The first term  in \eqref{K-21} simplifies to 
\begin{equation}\label{lastrecursionA}
\dfrac{[5]}{[2]}\kappa_{(a-1, b), (-1, 1)}. 
\end{equation}

We need to resolve the second and third terms on the right hand side of \eqref{K-21}. Both terms contain the following sub-diagram.
\begin{equation}\label{subdiag}
\begin{tikzpicture}[scale=.3, anchorbase]
	\draw[ultra thick, green] (1.5, -4) -- (1.5, -7);
	\draw[ultra thick, blue] (1.5,-4) arc (-90:0:1);
	\draw[ultra thick, blue] (0.5,-3) arc (-180:-90:1);
	\draw[ultra thick, blue] (2.5, -3) -- (2.5, 0);
	\draw[ultra thick, blue] (.5, -3) -- (.5, 0);
	\draw[ultra thick, black, fill= white] (-1.8,-2) ellipse (2.8 and 1.4) node{$a-1, b$};
	\draw[ultra thick, black, fill= white] (0,-7) ellipse (4 and 1) node{};
\end{tikzpicture}
\end{equation}
Expanding the $(a-1, b)$ clasp and using clasp orthogonality \eqref{clasporth} and neutral absorption \eqref{neutralabs}, we can rewrite \eqref{subdiag} as follows.
\begin{equation}\label{subdiagresolution}
\begin{tikzpicture}[scale=.3, anchorbase]
	\draw[ultra thick, green] (1.5, -4) -- (1.5, -7);
	\draw[ultra thick, blue] (1.5,-4) arc (-90:0:1);
	\draw[ultra thick, blue] (0.5,-3) arc (-180:-90:1);
	\draw[ultra thick, blue] (2.5, -3) -- (2.5, 0);
	\draw[ultra thick, blue] (.5, -3) -- (.5, 0);
	\draw[ultra thick, black, fill= white] (-1.8,-2) ellipse (2.8 and 1.4) node{$a-1, b$};
	\draw[ultra thick, black, fill= white] (0,-7) ellipse (4 and 1) node{};
\end{tikzpicture}
\quad = \quad
\begin{tikzpicture}[scale=.3, anchorbase]
	\draw[ultra thick, green] (1.5, -4) -- (1.5, -7);
	\draw[ultra thick, blue] (1.5,-4) arc (-90:0:1);
	\draw[ultra thick, blue] (0.5,-3) arc (-180:-90:1);
	\draw[ultra thick, blue] (2.5, -3) -- (2.5, 0);
	\draw[ultra thick, blue] (.5, -3) -- (.5, 0);
	\draw[ultra thick, black, fill= white] (-3,-2) ellipse (2.8 and 1.4) node{$a-2, b$};
	\draw[ultra thick, black, fill= white] (0,-7) ellipse (4 and 1) node{};
\end{tikzpicture}
\quad - \kappa_{(a-2, b), (-1, 1)}^{-1} \quad 
\begin{tikzpicture}[scale=.3, anchorbase]
	\draw[ultra thick, green] (.5, -4) -- (.5, -7);
	\draw[ultra thick, blue] (.5,-4) arc (-90:0:1);
	\draw[ultra thick, blue] (-0.5,-3) arc (-180:-90:1);
	\draw[ultra thick, blue] (1.5, -3) -- (1.5, 0);
	\draw[ultra thick, blue] (-.5, -3) -- (-.5, 0);
	\draw[ultra thick, blue] (2.5, -7) -- (2.5, 0);
	\draw[ultra thick, black, fill= white] (-2.8,-2) ellipse (2.8 and 1.4) node{$a-2, b$};
	\draw[ultra thick, black, fill= white] (-1,-7) ellipse (4 and 1) node{};
\end{tikzpicture}
\end{equation}
Then, from the relation \eqref{subdiagresolution}, we have
\begin{equation}
\begin{tikzpicture}[scale=.3, anchorbase]
	\draw[ultra thick, green] (1.5, -4) -- (1.5, -7);
	\draw[ultra thick, green] (1.5, 3) -- (1.5, 1);
	\draw[ultra thick, blue] (1.5,-4) arc (-90:0:1);
	\draw[ultra thick, blue] (0.5,-3) arc (-180:-90:1);
	\draw[ultra thick, blue] (2.5,0) arc (0:180:1);
	\draw[ultra thick, blue] (2.5, -3) -- (2.5, 0);
	\draw[ultra thick, blue] (.5, -3) -- (.5, 0);
	\draw[ultra thick, black, fill= white] (-1.8,-1.6) ellipse (2.8 and 1.4) node{$a-1, b$};
	\draw[ultra thick, black, fill= white] (0,-7) ellipse (4 and 1) node{};
\end{tikzpicture}
\quad = \quad
\begin{tikzpicture}[scale=.3, anchorbase]
	\draw[ultra thick, green] (1.5, -4) -- (1.5, -7);
	\draw[ultra thick, green] (1.5, 3) -- (1.5, 1);
	\draw[ultra thick, blue] (1.5,-4) arc (-90:0:1);
	\draw[ultra thick, blue] (0.5,-3) arc (-180:-90:1);
	\draw[ultra thick, blue] (2.5,0) arc (0:180:1);
	\draw[ultra thick, blue] (2.5, -3) -- (2.5, 0);
	\draw[ultra thick, blue] (.5, -3) -- (.5, 0);
	\draw[ultra thick, black, fill= white] (0,-7) ellipse (4 and 1) node{};
\end{tikzpicture}
\quad - \kappa_{(a-2, b), (-1, 1)}^{-1} \quad 
\begin{tikzpicture}[scale=.3, anchorbase]
	\draw[ultra thick, green] (2.5, -4) -- (2.5, 3);
	\draw[ultra thick, green] (-1, -7) -- (-1, -4);
	\draw[ultra thick, blue] (2.5, -7) -- (2.5, -4);
	\draw[ultra thick, blue] (-1, -4) -- (2.5, -4);
	\draw[ultra thick, blue] (-1, -4) -- (-1, -1);
	\draw[ultra thick, black, fill= white] (-2.8,-1.6) ellipse (2.8 and 1.4) node{$a-2, b$};
	\draw[ultra thick, black, fill= white] (-1,-7) ellipse (4 and 1) node{};
\end{tikzpicture}
\end{equation}
and
\begin{equation}
\begin{tikzpicture}[scale=.3, anchorbase]
	\draw[ultra thick, green] (1.5, -4) -- (1.5, -7);
	\draw[ultra thick, green] (.8, 3) -- (.8, .8);
	\draw[ultra thick, green] (2.2, 3) -- (2.2, .8);
	\draw[ultra thick, blue] (1.5,-4) arc (-90:0:1);
	\draw[ultra thick, blue] (0.5,-3) arc (-180:-90:1);
	\draw[ultra thick, blue] (2.5,0) arc (0:180:1);
	\draw[ultra thick, blue] (2.5, -3) -- (2.5, 0);
	\draw[ultra thick, blue] (.5, -3) -- (.5, 0);
	\draw[ultra thick, black, fill= white] (-1.8,-1.6) ellipse (2.8 and 1.4) node{$a-1, b$};
	\draw[ultra thick, black, fill= white] (0,-7) ellipse (4 and 1) node{};
\end{tikzpicture}
\quad = \quad
\begin{tikzpicture}[scale=.3, anchorbase]
	\draw[ultra thick, green] (1.5, -4) -- (1.5, -7);
	\draw[ultra thick, green] (.8, 3) -- (.8, .8);
	\draw[ultra thick, green] (2.2, 3) -- (2.2, .8);
	\draw[ultra thick, blue] (1.5,-4) arc (-90:0:1);
	\draw[ultra thick, blue] (0.5,-3) arc (-180:-90:1);
	\draw[ultra thick, blue] (2.5,0) arc (0:180:1);
	\draw[ultra thick, blue] (2.5, -3) -- (2.5, 0);
	\draw[ultra thick, blue] (.5, -3) -- (.5, 0);
	\draw[ultra thick, black, fill= white] (0,-7) ellipse (4 and 1) node{};
\end{tikzpicture}
\quad - \kappa_{(a-2, b), (-1, 1)}^{-1} \quad 
\begin{tikzpicture}[scale=.3, anchorbase]
	\draw[ultra thick, green] (2.5, -4) -- (2.5, 3);
	\draw[ultra thick, green] (-1, -7) -- (-1, -4);
	\draw[ultra thick, green] (2.5, -5) ..controls (3.5, -5) .. (3.5, 3);
	\draw[ultra thick, blue] (2.5, -7) -- (2.5, -4);
	\draw[ultra thick, blue] (-1, -4) -- (2.5, -4);
	\draw[ultra thick, blue] (-1, -4) -- (-1, -1);
	\draw[ultra thick, black, fill= white] (-2.8,-1.6) ellipse (2.8 and 1.4) node{$a-2, b$};
	\draw[ultra thick, black, fill= white] (-1,-7) ellipse (4 and 1) node{};
\end{tikzpicture}
\end{equation}
Applying these local relations to the second and third term on the right hand side of \eqref{K-21}, and simplifying diagrams using the defining relations of $\DD$, we obtain the next two equations. 
\begin{equation}\label{lastrecursionB}
\begin{tikzpicture}[scale=.25, anchorbase]
	\draw[ultra thick, green] (1.5, -4.5) -- (1.5, -7);
	\draw[ultra thick, green] (1.5, 4.5) -- (1.5, 7);
	\draw[ultra thick, green] (1.5, 1.5) -- (1.5, -1.5);
	\draw[ultra thick, blue] (1.5,-4.5) arc (270: 90: 1.5);
	\draw[ultra thick, blue] (1.5,-4.5) arc (-90: 90: 1.5);
	\draw[ultra thick, blue] (1.5,1.5) arc (270: 90: 1.5);
	\draw[ultra thick, blue] (1.5,1.5) arc (-90: 90: 1.5);
	\draw[ultra thick, black, fill= white] (0,7) ellipse (4 and 1) node{};
	\draw[ultra thick, black, fill= white] (-1,-3) ellipse (2.2 and .7) node{};
	\draw[ultra thick, black, fill= white] (-1,3) ellipse (2.2 and .7) node{};
	\draw[ultra thick, black, fill= white] (0,-7) ellipse (4 and 1) node{};
\end{tikzpicture}
\quad = - \left([2] + \kappa_{(a-2, b), (-1, 1)}^{-1}\right)\quad
\begin{tikzpicture}[scale=.25, anchorbase]
	\draw[ultra thick, green] (1.5, -2) -- (1.5, -7);
	\draw[ultra thick, green] (1.5, 2) -- (1.5, 7);
	\draw[ultra thick, blue] (1.5,-2) arc (270: 90: 2);
	\draw[ultra thick, blue] (1.5,-2) arc (-90: 90: 2);
	\draw[ultra thick, black, fill= white] (0,7) ellipse (4 and 1) node{};
	\draw[ultra thick, black, fill= white] (-1.5,0) ellipse (3 and 1) node{};
	\draw[ultra thick, black, fill= white] (0,-7) ellipse (4 and 1) node{};
\end{tikzpicture}
\end{equation}
\begin{equation}\label{lastrecursionC}
\begin{tikzpicture}[scale=.25, anchorbase]
	\draw[ultra thick, green] (1.5, -4.5) -- (1.5, -7);
	\draw[ultra thick, green] (1.5, 4.5) -- (1.5, 7);
	\draw[ultra thick, green] (.5, -1.9) ..controls (-.5, 0) .. (.5, 1.9);
	\draw[ultra thick, green] (2.5, -1.9) ..controls (3.5, 0) .. (2.5, 1.9);
	\draw[ultra thick, blue] (1.5,-4.5) arc (270: 90: 1.5);
	\draw[ultra thick, blue] (1.5,-4.5) arc (-90: 90: 1.5);
	\draw[ultra thick, blue] (1.5,1.5) arc (270: 90: 1.5);
	\draw[ultra thick, blue] (1.5,1.5) arc (-90: 90: 1.5);
	\draw[ultra thick, black, fill= white] (0,7) ellipse (4 and 1) node{};
	\draw[ultra thick, black, fill= white] (-1,-3) ellipse (2.2 and .7) node{};
	\draw[ultra thick, black, fill= white] (-1,3) ellipse (2.2 and .7) node{};
	\draw[ultra thick, black, fill= white] (-2,0) ellipse (2.2 and .7) node{};
	\draw[ultra thick, black, fill= white] (0,-7) ellipse (4 and 1) node{};
\end{tikzpicture}
\quad = \kappa_{(a-2, b), (-1, 1)}^{-2} \quad
\begin{tikzpicture}[scale=.25, anchorbase]
	\draw[ultra thick, green] (0, -4.5) -- (0, -7);
	\draw[ultra thick, green] (0, 4.5) -- (0, 7);
	\draw[ultra thick, green] (2, 4.5) -- (2, -4.5);
	\draw[ultra thick, green] (2, -5) ..controls (4.5, -7) and (4.5, 7) .. (2, 5);
	\draw[ultra thick, blue] (0, -4.5) -- (0, 4.5);
	\draw[ultra thick, blue] (2, -4.5) -- (0, -4.5);
	\draw[ultra thick, blue] (2, 4.5) -- (0, 4.5);
	\draw[ultra thick, blue] (2, -4.5) -- (2, -7);
	\draw[ultra thick, blue] (2, 4.5) -- (2, 7);
	\draw[ultra thick, black, fill= white] (0,7) ellipse (4 and 1) node{};
	\draw[ultra thick, black, fill= white] (0,0) ellipse (2.5 and 1) node{};
	\draw[ultra thick, black, fill= white] (0,-7) ellipse (4 and 1) node{};
\end{tikzpicture}
\end{equation}

Note that the diagram on the right hand side of \eqref{lastrecursionB} is $K_{(a-1, b), (-1, 1)}$. Moreover, after applying neutral absorption \eqref{neutralabs} the diagram on the right hand side of \eqref{lastrecursionC} is $K_{(a-2, b+1), (0, 0)}$. Therefore, we can use \eqref{lastrecursionA}, \eqref{lastrecursionB}, and \eqref{lastrecursionC} to rewrite \eqref{K-21}, then deduce that 
\begin{equation}
\begin{split}
\kappa_{(a,b), (-2,1)} = \dfrac{[5]}{[2]} \kappa_{(a-1, b), (-1, 1)} &- (-[2]- \kappa_{(a-2,b), (-1,1)}^{-1})\dfrac{\kappa_{(a-1,b),(-1, 1)}}{\kappa_{(a-1, b),(-1,0)}} \\
&- \dfrac{\kappa_{(a-2, b+1), (0,0)}}{\kappa_{(a-2, b), (-1, 1)}^2 \kappa_{(a-1, b), (-1, 1)}}.
\end{split}
\end{equation}
\end{proof}

\subsection{Solving the recursion}
\label{subsec-solve}

\begin{prop}\label{formulas}
The recursive relations in Proposition \eqref{recursionprop} together with the initial conditions in Equations \eqref{ic1}, \eqref{ic2}, \eqref{ic3}, \eqref{ic4}, and \eqref{ic5}
are uniquely solved by 
\begin{equation}
\kappa_{(a,b), (1, 0)}= 1
\end{equation}
\begin{equation}
\kappa_{(a, b), (0,1)}= 1
\end{equation}
\begin{equation}\label{rs1}
\kappa_{(a, b), (-1, 1)} = - \dfrac{[a+1]}{[a]}
\end{equation}
\begin{equation}\label{rs2}
\kappa_{(a,b), (2, -1)} = - \dfrac{[2b+2]}{[2b]}. 
\end{equation}
\begin{equation}\label{rs3}
\kappa_{(a,b), (0,0)} = \dfrac{[a+2][a+2b+4]}{[2][a][a+2b+2]}. 
\end{equation}
\begin{equation}\label{rs4}
\kappa_{(a,b), (1, -1)} = \dfrac{[a+2b+3][2b+2]}{[a+ 2b+2][2b]}.
\end{equation}
\begin{equation}\label{rs5}
\kappa_{(a,b), (-2, 1)} = -\dfrac{[a+1][2a+ 2b+4]}{[a-1][2a+ 2b+ 2]}. 
\end{equation}
\begin{equation}\label{rs6}
\kappa_{(a,b), (-1, 0)} = - \dfrac{[2a+ 2b+ 4][a+ 2b+ 3][a+ 1]}{[2a+ 2b+ 2][a+ 2b+ 2][a]}. 
\end{equation}
\begin{equation}\label{rs7}
\kappa_{(a,b), (0, -1)} = \dfrac{[2a+2b+4][a+2b+3][2b+2]}{[2a+2b+2][a+2b+1][2b]}
\end{equation}
\end{prop}

\begin{proof}
There is a recursive relation for each non-dominant weight in a fundamental representation. We say that the right hand side of a relation \emph{involves} the weight $\mu$ if $\kappa_{?, \mu}$ appears in the right hand side of the recursion. 

That relation \eqref{rr1} (with the specified initial conditions) is solved by \eqref{rs1} is easily seen to be equivalent to showing that
\begin{equation}\label{easy1}
-[a+1] = -[2][a] + [a-1].
\end{equation}
This is a well known identity for quantum numbers, but we will describe a different way to derive \eqref{easy1}. First, we multiply equation \eqref{easy1} by $(q-q^{-1})$, resulting in
\begin{equation}\label{easyq}
-\left(q^{a+1}- q^{-(a+1)}\right) = -(q+q^{-1})\left(q^a- q^{-a}\right) + \left(q^{(a-1)} - q^{(-(a-1)}\right). 
\end{equation}
Second, we temporarily replace $q^a$ with the variable $A$, so \eqref{easyq} becomes the following.
\begin{equation}\label{easytocheck}
Aq- A^{-1}q^{-1} = -(q+q^{-1})(A - A^{-1}) + (Aq^{-1}- A^{-1}q)
\end{equation}
Equation \eqref{easytocheck} is easily seen to be true in $\mathbb{Z}[A^{\pm 1} , q^{\pm 1}]$, then specializing $A$ to $q^a$ we find that \eqref{easyq} holds as well. To see that relation \eqref{rr2} is solved by \eqref{rs1} is similar, and we leave it as an exercise. 

By using
\begin{equation}
\kappa_{(a, b), (-1, 1)} = -\dfrac{[a+1]}{[a]} \ \ \ \ \ \ \text{and} \ \ \ \ \ \kappa_{(a, b), (2 -1)} = -\dfrac{[2b+2]}{[2b]},
\end{equation}
we can simplify \eqref{rr3} and \eqref{rr4} as follows.
\begin{equation}
\kappa_{(a, b), (0, 0)} = \dfrac{[5]}{[2]} - \dfrac{[a-1]}{[a]} \dfrac{[2b+4]}{[2b+2]} - \kappa_{(a-1, b), (1, -1)}^{-1}
\end{equation}
\begin{equation}
\kappa_{(a, b), (1, -1)} = \dfrac{[5]}{[2]} -\dfrac{[2b-2]}{[2b]} \dfrac{[a+3]}{[a+2]}- \dfrac{1}{[2]^2}\kappa_{(a, b-1), (0, 0)}^{-1}
\end{equation}
Then by induction our claim that \eqref{rr3} is solved by \eqref{rs3} and \eqref{rr4} is solved by \eqref{rs4} follows from verifying the following two equalities.
\begin{equation}\label{desiredequality1}
\dfrac{[a+2][a+2b+4]}{[2][a][a+2b+2]}= \dfrac{[5]}{[2]} -  \dfrac{[a-1][2b+4]}{[a][2b+2]} - \dfrac{[a+2b+1][2b]}{[a+ 2b+2][2b+2]}
\end{equation}
\begin{equation}\label{desiredequality2}
\dfrac{[a+2b+3][2b+2]}{[a+ 2b+2][2b]}=  \dfrac{[5]}{[2]} -\dfrac{[a+3][2b-2]}{[a+2][2b]}- \dfrac{1}{[2]^2}\dfrac{[2][a][a+2b]}{[a+2][a+2b+2]}
\end{equation}

We focus on the quantum number calculation needed to verify the first of these two equalities. After clearing denominators the desired equality \eqref{desiredequality1} will follow from the following identity. 
\begin{equation}\label{desiredequality1a}
\begin{split}
[a+2][2b+2][a+2b+4] &= 
[5][a][2b+2][a+2b+2] \\ &- [2][a-1][2b+4][a+2b+2] - [2][a][2b][a+2b+1].
\end{split}
\end{equation}
Multiplying through by $(q-q^{-1})^3$ and replacing $q^a$ with $A$ and $q^b$ with $B$, we find the desired quantum number identity is a consequence of the following identity in $\mathbb{Z}[A^{\pm 1}, B^{\pm 1}, q^{\pm 1}]$. 
\begin{equation}\label{easytocheck2}
\begin{split}
(Aq^2- A^{-1}q^{-2})(B^2q^2- B^{-2}q^{-2})&(AB^2q^4- A^{-1}B^{-2}q^{-4}) = \\
(q^4+ q^2+ 1+ q^{-2} + q^{-4})(A- A^{-1})&(B^2q^2- B^{-2}q^{-2})(AB^2q^2-A^{-1}B^{-2}q^{-2}) \\
-(q+q^{-1})(Aq^{-1}- A^{-1}q)&(B^2q^4- B^{-2}q^{-4})(AB^2q^2- A^{-1}B^{-2}q^{-2}) \\
- (q+q^{-1})(A-A^{-1})&(B^2- B^{-2})(AB^2q- A^{-1}B^{-2}q^{-1})
\end{split}
\end{equation}
We leave the details of checking \eqref{easytocheck2} by hand as an exercise for the reader. Then replacing $A$ with $q^a$ and $B$ with $q^b$ we may deduce the equality \eqref{desiredequality1a}. 

The calculations needed to verify \eqref{desiredequality2} are omitted, as are the rest of the details of the quantum number calculations. We simply outline the remainder of the proof below. 

Once the first four relations are solved, we can simplify the fifth relation so the simplified recursion only involves the weight $(-1, 0)$. By using induction we reduce proving the recursion relation \eqref{rr5} is solved by \eqref{rs5} to a quantum number calculation. The sixth recursion relation only involves the previous five weights, so we can use these solutions to simplify the right hand side of \eqref{rr6}. A quantum number calculation will verify that the right hand side is in fact equal to \eqref{rs6}. After using the first six solutions to simplify the last recursion, \eqref{rr7} only involves the weight $(0, -1)$ and so can be solved by induction and a quantum number calculation. 
\end{proof}

\begin{rmk}
The point of writing \eqref{easy1} in the form \eqref{easytocheck} and \eqref{desiredequality1} in the form \eqref{easytocheck2} is that it makes it possible (and in fact easy) to have a computer verify the desired quantum number identities \cite{bodishc2code}. 
\end{rmk}

\subsection{Relation to Elias's clasp conjecture}
\label{subsec-claspconj}

In the following, we will reinterpret Elias's type $A$ clasp conjecture \cite{elias2015light} in type $C_2$. We then discuss how we expect Elias's conjecture generalizes to a type independent statement.

Recall that the Weyl group for the $C_2$ root system, which we denote simply by $W$, acts on the weight lattice $X$ by
\begin{equation}
s(\varpi_1)= - \varpi_1 + \varpi_2 \ \ \ \ \ \text{and} \ \ \ \ \ s(\varpi_2) = \varpi_2
\end{equation}
 while 
 \begin{equation}
 t(\varpi_1) = \varpi_1 \ \ \ \ \ \text{and} \ \ \ \ \ t(\varpi_2) = 2\varpi_1 - \varpi_2.
 \end{equation}
 
For a weight $\mu$, we will denote by $d_{\mu}$ the minimal length element in $W$ which when takes $\mu$ to a dominant weight. Thus, 
\begin{equation}
d_{\varpi_1} = 1 = d_{\varpi_2}
\end{equation}
and
\begin{equation}
d_{- \varpi_1 +\varpi_2} = s, d_{2\varpi_1 - \varpi_2} = t, d_{\varpi_1 - \varpi_2} = st, d_{-2\varpi_1 + \varpi_2} = ts, d_{- \varpi_1} = sts, \text{and} \ d_{- \varpi_2} = tst.
\end{equation}
We then define the set $\Phi_{\mu} = \lbrace \alpha\in \Phi_+ \ | \ d_{\mu}(\alpha)\in \Phi_-\rbrace$. Thus, 
\begin{equation}
\Phi_{- \varpi_1 + \varpi_2} = \lbrace \alpha_s\rbrace 
\end{equation}
\begin{equation}
\Phi_{2\varpi_1 - \varpi_2} = \lbrace \alpha_t\rbrace 
\end{equation}
\begin{equation}
\Phi_{\varpi_1 - \varpi_2} = \lbrace \alpha_t, t(\alpha_s)\rbrace 
\end{equation}
\begin{equation}
\Phi_{-2\varpi_1 + \varpi_2} = \lbrace \alpha_s, s(\alpha_t) \rbrace
\end{equation}
\begin{equation}
\Phi_{- \varpi_1} = \lbrace \alpha_s, s(\alpha_t), st(\alpha_s)\rbrace 
\end{equation}
\begin{equation}
\Phi_{-\varpi_2} = \lbrace \alpha_t, t(\alpha_s), ts(\alpha_t)\rbrace 
\end{equation}

Let $(-,-)$ be the standard inner product on $X$ so the $\epsilon_i$ are an orthonormal basis. Recall that $\alpha^{\vee} = 2\alpha/(\alpha, \alpha)$, and that $\rho$ is the sum of the fundamental weights. We define $q_{\alpha}= q$ when $\alpha$ is a short root and $q_{\alpha}= q^2$ when $\alpha$ is a long root.

\begin{cor}
In type $C_2$, if $\mu$ is an (extremal) weight in a fundamental representation and $\lambda\in X_+$, then 
\begin{equation}
\kappa_{\lambda, \mu} = \pm\prod_{\alpha\in \Phi_{\mu}} \dfrac{[(\alpha^{\vee}, \lambda + \rho)]_{q_{\alpha}}}{[(\alpha^{\vee}, \lambda + \mu + \rho)]_{q_{\alpha}}}
\end{equation}
\end{cor}

\begin{proof} 
Using the formula $\dfrac{[2n]_v}{[2]_v} = [n]_{v^2}$ it is an easy exercise to use \eqref{formulas} and our description of $\Phi_{\mu}$ to check the corollary.  
\end{proof}

It is natural to expect Elias's clasp conjecture to generalize as follows. Let $\Phi$ be an irreducible root system with associated simple Lie algebra $\mathfrak{g}$. Let $(-, -)$ be the $W$ invariant bilinear form on $\Phi$ so that $(\alpha, \alpha)= 2$ for all short roots $\alpha\in \Phi$. Fix a fundamental weight $\varpi$ and a weight $\mu\in \wt V(\varpi)$ which is in the Weyl group orbit of $\varpi$. Let $d_{\mu}$ be the minimal length element in the Weyl group such that $d_{\mu}(\mu) = \varpi$. Set 
\[
\Phi_{+}(\mu)= \lbrace \alpha \ | \ d_{\mu}(\alpha)\in \Phi_{-}\rbrace.  
\]

\begin{conj}\label{claspconjecture}
There is an elementary light ladder map $L_{\mu}$, which is a morphism of $U_q(\mathfrak{g})$ modules, and for each dominant weight $\lambda$ a map $E_{\lambda, \mu}$ (which may be zero) which is a composition of the clasps $C_{\lambda}$ and $C_{\lambda+ \mu}$ with $L_{\mu}$ as in Equation \eqref{Emap}. Moreover, there is a duality $\mathbb{D}$ on $\Fund(\mathfrak{g})$ which, interpreted in the graphical calculus, is flipping a diagram upside down. Finally, we expect that
\begin{equation}\label{generalLIF}
E_{\lambda, \mu}\circ \mathbb{D}(E_{\lambda, \mu}) = \prod_{\alpha\in \Phi_+(\mu)}\dfrac{[(\alpha^{\vee}, \lambda+ \rho)]_{q_{\alpha}}}{[(\alpha^{\vee}, \lambda+ \rho + \mu)]_{q_{\alpha}}} C_{\lambda},
\end{equation}
where $q_{\alpha}= q^{(\alpha, \alpha)/2}$. 
\end{conj}

We can already conjecture the general form of one of the recursive relations satisfied by the local intersection forms. The local intersection form calculations in type $C_2$ show that every weight in $V(\varpi)$ appears in this recursion for the local intersection form of $\kappa_{\lambda, - \varpi}$. 
\begin{conj}\label{lowestweightrecursion}
\begin{equation}
\kappa_{\lambda, - \varpi} = \dim_q V(\varpi) - \sum_{\mu\in \wt V(\varpi)} \dfrac{\kappa_{\lambda- \varpi + \mu, - \mu}}{\kappa_{\lambda - \varpi, \mu}}.
\end{equation}
\end{conj}
One might hope to prove Conjecture \eqref{claspconjecture} by finding a combinatorial description of the recursive formula's themselves, then proving these recursions are both given by calculations with webs and  solved uniquely by Equation \eqref{generalLIF}. 

\begin{rmk}
The conjecture in type $A$ only deals with $\kappa_{\lambda, \mu}$ when $\mu$ is in the Weyl group orbit of a dominant fundamental weight. In type $C_2$, we cannot currently explain the local intersection form for the weight $(0,0)$ in a way that suggests any generalization. However, we do expect there is a general formula which, for any simple Lie algebra $\mathfrak{g}$ and any fundamental weight $\varpi$, computes $\kappa_{\lambda, \mu}$ for all $\mu \in \wt V(\varpi)$ in terms of the root system $\Phi$. 
\end{rmk}

\bibliographystyle{plain}

\bibliography{mastercopy}



\end{document}